\documentclass[11pt]{article}
\usepackage{amssymb}
\usepackage{amsmath}
\usepackage[all]{xy}

\title{\bf Wakamatsu Tilting Modules, $U$-Dominant Dimension and $k$-Gorenstein Modules
}
\author{Zhaoyong Huang\thanks{\small \it E-mail address: huangzy@nju.edu.cn}\\
{\small \it Department of Mathematics, Nanjing University,}\\
{\small \it Nanjing 210093, People's Republic of China}\\
\\
{\footnotesize \it Dedicated to Professor Edgar E. Enochs on his
72nd birthday}}
\usepackage{amssymb}
\date{}
\begin{document}
\baselineskip=18pt \maketitle

\begin{abstract}
Let $\Lambda$ and $\Gamma$ be left and right noetherian rings and
$_{\Lambda}U$ a Wakamatsu tilting module with $\Gamma ={\rm
End}(_{\Lambda}T)$. We introduce a new definition of $U$-dominant
dimensions and show that the $U$-dominant dimensions of
$_{\Lambda}U$ and $U _{\Gamma}$ are identical. We characterize
$k$-Gorenstein modules in terms of homological dimensions and the
property of double homological functors preserving monomorphisms.
We also study a generalization of $k$-Gorenstein modules, and
characterize it in terms of some similar properties of
$k$-Gorenstein modules.
\end{abstract}

\vspace{1cm}

\centerline{\bf 1. Introduction and main results}

\vspace{0.2cm}

Let $\Lambda$ be a ring. We use Mod $\Lambda$ (resp. Mod $\Lambda
^{op}$) to denote the category of left (resp. right)
$\Lambda$-modules, and use mod $\Lambda$ (resp. mod $\Lambda
^{op}$) to denote the category of finitely generated left
$\Lambda$-modules (resp. right $\Lambda$-modules).

\vspace{0.2cm}

{\bf Definition 1.1}$^{[7]}$ For a module $M$ in mod $\Lambda$
(resp. mod $\Lambda ^{op}$) and a positive integer $k$, $M$ is
said to have dominant dimension at least $k$, written as
dom.dim$(_{\Lambda}M)$ (resp. dom.dim$(M_\Lambda)) \geq k$, if
each of the first $k$ terms in a minimal injective resolution of
$M$ is $\Lambda$-flat (resp. $\Lambda ^{op}$-flat).

\vspace{0.2cm}

For a module $T$ in Mod $\Lambda$ (resp. Mod $\Lambda ^{op}$), we
use add-lim$_{\Lambda}T$ (resp. add-lim$T_{\Lambda}$) to denote
the subcategory of Mod $\Lambda$ (resp. Mod $\Lambda ^{op}$)
consisting of all modules isomorphic to direct summands of a
direct limit of a family modules in which each is a finite direct
sum of copies of $_{\Lambda}T$ (resp. $T_{\Lambda}$). We now
introduce a definition of $U$-dominant dimension as follows.

\vspace{0.2cm}

{\bf Definition 1.2} Let $U$ be in Mod $\Lambda$ (resp. Mod
$\Lambda ^{op}$) and $k$ a positive integer. For a module $M$ in
Mod $\Lambda$ (resp. Mod $\Lambda ^{op}$), $M$ is said to have
$U$-dominant dimension at least $k$, written as
$U$-dom.dim$(_{\Lambda}M)$ (resp. $U$-dom.dim$(M_{\Lambda}))\geq
k$, if each of the first $k$ terms in a minimal injective
resolution of $M$ can be embedded into a direct limit of a family
of modules in which each is a finite direct sum of copies of
$_{\Lambda}U$ (resp. $U _{\Lambda}$), that is, each of these terms
is in add-lim$_{\Lambda}U$ (resp. add-lim$U_{\Lambda}$).

\vspace{0.2cm}

{\it Remark.} Notice that a module (not necessarily finitely
generated) is flat if and only if it is a direct limit of a family
of finitely generated free modules (see [15]). So, if putting
$_{\Lambda}U={_{\Lambda}\Lambda}$ (resp. $U_{\Lambda}=\Lambda
_{\Lambda}$), then the above definition of $U$-dominant dimension
coincides with that of the usual dominant dimension for any ring
$\Lambda$.

\vspace{0.2cm}

Tachikawa in [19] showed that if $\Lambda$ is a left and right
artinian ring then the dominant dimensions of $_{\Lambda}\Lambda$
and $\Lambda _{\Lambda}$ are identical. Hoshino in [7] further
showed that this result also holds for left and right noetherian
rings. Colby and Fuller in [5] gave some equivalent conditions of
dom.dim$(_{\Lambda}\Lambda) \geq 1$ (or 2) in terms of the
properties of double dual functors (with respect to
${_{\Lambda}\Lambda} _{\Lambda}$). These results motivate our
interests in establishing the identity of $U$-dominant dimensions
of $_{\Lambda}U$ and $U _{\Gamma}$ (where $\Gamma ={\rm
End}(_{\Lambda}U)$) and characterizing the properties of modules
with a given $U$-dominant dimension.

Let $T$ be a module in mod $\Lambda$. For a module $A \in$ mod
$\Lambda$ and a non-negative integer $n$, we say that the grade of
$A$ with respect to $_{\Lambda}T$, written as grade$_TA$, is at
least $n$ if Ext$_{\Lambda}^i(A, T)=0$ for any $0 \leq i < n$. We
say that the strong grade of $A$ with respect to $T$, written as
s.grade$_TA$, is at least $n$ if grade$_TB \geq n$ for all
submodules $B$ of $A$. The notion of the (strong) grade of modules
with respect to a given module in mod $\Lambda ^{op}$ is defined
dually.

The following is one of main results in this paper.

\vspace{0.2cm}

{\bf Theorem I} {\it Let} $\Lambda$ {\it and} $\Gamma$ {\it be
left and right noetherian rings and} $_{\Lambda}U$ {\it a
Wakamatsu tilting module with} $\Gamma ={\rm End}(_{\Lambda}U)$.
{\it For a positive integer} $k$, {\it the following statements
are equivalent}.

(1) $U$-dom.dim$(_{\Lambda}U) \geq k$.

(2) s.grade$_U$Ext$_{\Lambda}^1(M, U) \geq k$ {\it for any} $M
\in$mod $\Lambda$.

(3) Hom$_{\Lambda}(U, E_i)$ {\it is} $\Gamma$-{\it flat, where}
$E_i$ {\it is the} $(i+1)$-st {\it term in a minimal injective
resolution of} $_{\Lambda}U$, {\it for any} $0 \leq i \leq k-1$.

(1)$^{op}$ $U$-dom.dim$(U_{\Gamma}) \geq k$.

(2)$^{op}$ s.grade$_U$Ext$_{\Gamma}^1(N, U) \geq k$ {\it for any}
$N\in$mod $\Gamma ^{op}$.

(3)$^{op}$ Hom$_{\Gamma}(U, E'_i)$ {\it is} $\Lambda ^{op}$-{\it
flat, where} $E'_i$ {\it is the} $(i+1)$-st {\it term in a minimal
injective resolution of} $U_{\Gamma}$, {\it for any} $0 \leq i
\leq k-1$.

\vspace{0.2cm}

Kato in [14] gave a definition of $U$-dominant dimension as
follows, which is different from that of Definition 1.2. For a
module $M$ in mod $\Lambda$ (resp. mod $\Lambda ^{op}$), $M$ is
said to have $U$-dominant dimension at least $k$, written as
$U$-dom.dim$(_{\Lambda}M)$ (resp. $U$-dom.dim$(M_{\Lambda}))\geq
k$, if each of the first $k$ terms in a minimal injective
resolution of $M$ is cogenerated by $_{\Lambda}U$ (resp. $U
_{\Lambda}$), that is, each of these terms can be embedded into a
direct product of copies of $_{\Lambda}U$ (resp. $U _{\Lambda}$).
If we adopt the definition of $U$-dominant dimension given by
Kato, then in Theorem I the equivalence of (2), (3), (2)$^{op}$
and (3)$^{op}$ and that (1) implies (3) also hold. However, that
(3) does not imply (1) in general. For example, consider Wakamatsu
tilting module $_{\mathbb{Z}}\mathbb{Z}$ and its injective
envelope $_{\mathbb{Z}}\mathbb{Q}$, where $\mathbb{Z}$ is the ring
of integers and $\mathbb{Q}$ is the field of rational numbers.
Then the module $_{\mathbb{Z}}\mathbb{Q}$ is flat, but it can not
be embedded into any direct product of copies of
$_{\mathbb{Z}}\mathbb{Z}$ since Hom$_{\mathbb{Z}}(\mathbb{Q},
\mathbb{Z})=0$.

\vspace{0.2cm}

{\bf Corollary 1.3} {\it Let} $\Lambda$ {\it and} $\Gamma$ {\it be
left and right noetherian rings and} $_{\Lambda}U$ {\it a
Wakamatsu tilting module with} $\Gamma ={\rm End}(_{\Lambda}T)$.
{\it Then} $U$-dom.dim$(_{\Lambda}U)=U$-dom.dim$(U _{\Gamma})$.

\vspace{0.2cm}

{\it Remark.} We do not know whether the conclusion in Corollary
1.3 holds for Kato's $U$-dominant dimension. The answer is
positive when $\Lambda$ and $\Gamma$ are artinian algebras (see
[11] Theorem 1.3).

\vspace{0.2cm}

Putting $_{\Lambda}U_{\Gamma}={_{\Lambda}\Lambda _{\Lambda}}$, we
immediately get the following result, which is due to Hoshino (see
[7] Theorem).

\vspace{0.2cm}

{\bf Corollary 1.4} {\it For a left and right noetherian ring}
$\Lambda$, dom.dim$(_{\Lambda}\Lambda)=$dom.dim$(\Lambda
_{\Lambda})$.

\vspace{0.2cm}

{\bf Definition 1.5}$^{[12]}$ For a non-negative integer $k$, a
module $U \in$ mod $\Lambda$ with $\Gamma ={\rm End}(_{\Lambda}U)$
is called $k$-Gorenstein if s.grade$_{U}{\rm Ext}^{i}_{\Gamma} (N,
U)\geq i$ for any $N \in$mod $\Gamma ^{op}$ and $1 \leq i \leq k$.
Dually, we may define the notion of $k$-Gorenstein modules in mod
$\Gamma ^{op}$.

\vspace{0.2cm}

We introduce a new homological dimension of modules as follows.

\vspace{0.2cm}

{\bf Definition 1.6} Let $\Lambda$ be a ring and $T$ in Mod
$\Lambda$. For a module $A$ in Mod $\Lambda$, if there is an exact
sequence $\cdots \to T_{n} \to \cdots \to T_{1} \to T_{0} \to A
\to 0$ in Mod $\Lambda$ with each $T_{i}\in$add-lim$_{\Lambda}T$
for any $i \geq 0$, then we define
$T$-lim.dim$_{\Lambda}(A)=$inf$\{ n|$ there is an exact sequence
$0 \to T_{n} \to \cdots \to T_{1} \to T_{0} \to A \to 0$ in Mod
$\Lambda$ with each $T_{i}\in$add-lim$_{\Lambda}T$ for any $0 \leq
i \leq n \}$. We set $T$-lim.dim$_{\Lambda}(A)$ infinity if no
such an integer exists. For $\Lambda^{op}$-modules, we may define
such a dimension dually.

\vspace{0.2cm}

{\it Remark.} Putting $_{\Lambda}T={_{\Lambda}\Lambda}$ (resp.
$T_{\Lambda}=\Lambda_{\Lambda}$), the dimension defined as above
is just the flat dimension of modules.

\vspace{0.2cm}

In [21], Wakamatsu showed that the notion of $k$-Gorenstein
modules is left-right symmetric. We give here some other
characterizations of $k$-Gorenstein modules. The following is
another main result in this paper.

\vspace{0.2cm}

{\bf Theorem II} {\it Let} $\Lambda$ {\it and} $\Gamma$ {\it be
left and right noetherian rings and} $_{\Lambda}U$ {\it a
Wakamatsu tilting module with} $\Gamma ={\rm End}(_{\Lambda}T)$.
{\it Then, for a positive integer} $k$, {\it the following
statements are equivalent}.

(1) $_{\Lambda}U$ {\it is} $k$-{\it Gorenstein}.

(2) s.grade$_U$Ext$_{\Lambda}^i(M, U) \geq i$ {\it for any} $M
\in$mod $\Lambda$ {\it and} $1 \leq i \leq k$.

(3) $U$-lim.dim$_{\Lambda}(E_i) \leq i$ {\it for any} $0 \leq i
\leq k-1$.

(4) {\it l.}fd$_{\Gamma}($Hom$_{\Lambda}(U, E_i)) \leq i$ {\it for
any} $0 \leq i \leq k-1$, {\it where l.}fd {\it denotes the left
flat dimension and} $E_i$ {\it is the} $(i+1)$-st {\it term in a
minimal injective resolution of} $_{\Lambda}U$, {\it for any} $0
\leq i \leq k-1$.

(5) ${\rm Ext}_{\Gamma}^i({\rm Ext}_{\Lambda}^i(\ , U), U)$ {\it
preserves monomorphisms in} mod $\Lambda$ {\it for any} $0 \leq i
\leq k-1$.

(1)$^{op}$ $U_{\Gamma}$ {\it is} $k$-{\it Gorenstein}.

(2)$^{op}$ s.grade$_U$Ext$_{\Gamma}^i(N, U) \geq i$ {\it for any}
$N \in$mod $\Gamma ^{op}$ {\it and} $1 \leq i \leq k$.

(3)$^{op}$ $U$-lim.dim$_{\Gamma}(E'_i) \leq i$ {\it for any} $0
\leq i \leq k-1$.

(4)$^{op}$ {\it r.}fd$_{\Lambda}($Hom$_{\Gamma}(U, E'_i)) \leq i$
{\it for any} $0 \leq i \leq k-1$, {\it where r.}fd {\it denotes
the right flat dimension and} $E'_i$ {\it is the} $(i+1)$-st {\it
term in a minimal injective resolution of} $U_{\Gamma}$ {\it for
any} $0 \leq i \leq k-1$.

(5)$^{op}$ ${\rm Ext}_{\Lambda}^i({\rm Ext}_{\Gamma}^i(\ , U), U)$
{\it preserves monomorphisms in} mod $\Gamma ^{op}$, {\it for any}
$0 \leq i \leq k-1$.

\vspace{0.2cm}

Let $\Lambda$ and $\Gamma$ be left and right noetherian rings and
$_{\Lambda}U$ a Wakamatsu tilting module with $\Gamma ={\rm
End}(_{\Lambda}U)$. By Theorems I and II, if $U$ has $U$-dominant
dimension at least $k$, then it is $k$-Gorenstein.

Recall that a left and right noetherian ring $\Lambda$ is called
$k$-Gorenstein if the flat dimension of the $i$-th term in a
minimal injective resolution of $_{\Lambda}\Lambda$ is at most
$i-1$ for any $1 \leq i \leq k$. Auslander showed in [6] Theorem
3.7 that the notion of $k$-Gorenstein rings is left-right
symmetric. Following Definition 1.6 and [6] Theorem 3.7, a left
and right noetherian ring $\Lambda$ is $k$-Gorenstein if it is
$k$-Gorenstein as a $\Lambda$-module. So, by Theorem II, we have
the following corollary, which develops this Auslander's result.

\vspace{0.2cm}

{\bf Corollary 1.7} {\it Let} $\Lambda$ {\it and} $\Gamma$ {\it be
left and right noetherian rings. Then, for a positive integer}
$k$, {\it the following statements are equivalent}.

(1) $\Lambda$ {\it is} $k$-{\it Gorenstein}.

(2) s.grade$_{\Lambda}$Ext$_{\Lambda}^i(M, \Lambda) \geq i$ {\it
for any} $M \in$mod $\Lambda$ {\it and} $1 \leq i \leq k$.

(3) {\it The flat dimension of the} $i$-th {\it term in a minimal
injective resolution of} $_{\Lambda}\Lambda$ {\it is at most}
$i-1$ {\it for any} $1 \leq i \leq k$.

(4) ${\rm Ext}_{\Lambda}^i({\rm Ext}_{\Lambda}^i(\ , \Lambda),
\Lambda)$ {\it preserves monomorphisms in} mod $\Lambda$ {\it for
any} $0 \leq i \leq k-1$.

(2)$^{op}$ s.grade$_{\Lambda}$Ext$_{\Lambda}^i(N, \Lambda) \geq i$
{\it for any} $N \in$mod $\Lambda ^{op}$ {\it and} $1 \leq i \leq
k$.

(3)$^{op}$ {\it The flat dimension of the} $i$-th {\it term in a
minimal injective resolution of} $\Lambda _{\Lambda}$ {\it is at
most} $i-1$ {\it for any} $1 \leq i \leq k$.

(4)$^{op}$ ${\rm Ext}_{\Lambda}^i({\rm Ext}_{\Lambda}^i(\ ,
\Lambda), \Lambda)$ {\it preserves monomorphisms in} mod $\Lambda
^{op}$ {\it for any} $0 \leq i \leq k-1$.

\vspace{0.2cm}

The paper is organized as follows. In Section 2, we give some
properties of Wakamatsu tilting modules. For example, let
$\Lambda$ and $\Gamma$ be left and right noetherian rings and
$_{\Lambda}U$ a Wakamatsu tilting module with $\Gamma ={\rm
End}(_{\Lambda}U)$. If $_{\Lambda}U$ is $k$-Gorenstein for all
$k$, then the left and right injective dimensions of
$_{\Lambda}U_{\Gamma}$ are identical provided that both of them
are finite. We shall prove our main results in Section 3. As
applications of the results obtained in Section 3, we characterize
in Section 4 $U$-dominant dimension of $U$ at least one and two in
terms of the properties of ${\rm Hom}({\rm Hom}(\ , U), U)$
preserving monomorphisms and being left exact, respectively.
Motivated by the work of Auslander and Reiten in [3], we study in
Section 5 a generalization of $k$-Gorenstein modules, which is
however not left-right symmetric. We characterize this
generalization in terms of some properties similar to that of
$k$-Gorenstein modules. At the end of this section, we generalize
the result of Wakamatsu on the symmetry of $k$-Gorenstein modules.

\vspace{0.5cm}

\centerline{\bf 2. Wakamatsu tilting modules}

\vspace{0.2cm}

In this section, we give some properties of Wakamatsu tilting
modules with finite homological dimensions.

\vspace{0.2cm}

{\bf Definition 2.1} Let $\Lambda$ be a ring. A module
$_{\Lambda}U$ in mod $\Lambda$ is called a Wakamatsu tilting
module if $_{\Lambda}U$ is self-orthogonal (that is, ${\rm
Ext}_{\Lambda}^{i}(_{\Lambda}U , {_{\Lambda}U})=0$ for any $i \geq
1$), and possessing an exact sequence:
$$0 \to {_{\Lambda}\Lambda} \to U_0 \to U_1 \to \cdots \to U_i \to \cdots$$
such that: (1) all term $U_i$ are direct summands of finite direct
sums of copies of $_{\Lambda}U$, that is, $U_i
\in$add$_{\Lambda}U$, and (2) after applying the functor
Hom$_{\Lambda}(\ , U)$ the sequence is still exact. The definition
of Wakamatsu tilting modules in mod $\Lambda ^{op}$ is given
dually (see [20] or [21]).

\vspace{0.2cm}

Let $\Lambda$ and $\Gamma$ be rings. Recall that a bimodule
$_{\Lambda}U _{\Gamma}$ is called a faithfully balanced
selforthogonal bimodule if it satisfies the following conditions:

(1) $_{\Lambda}U \in$mod $\Lambda$ and $U_{\Gamma} \in$mod $\Gamma
^{op}$.

(2) The natural maps $\Lambda \rightarrow {\rm End}(U_{\Gamma})$
and $\Gamma \rightarrow {\rm End}(_{\Lambda}U)^{op}$ are
isomorphisms.

(3) ${\rm Ext}_{\Lambda}^{i}(_{\Lambda}U , {_{\Lambda}U})=0$ and
${\rm Ext}_{\Gamma}^{i}(U_{\Gamma}, U_{\Gamma})=0$ for any $i \geq
1$.

\vspace{0.2cm}

The following result is [21] Corollary 3.2.

\vspace{0.2cm}

{\bf Proposition 2.2} {\it Let} $\Lambda$ {\it be a left
noetherian ring and} $\Gamma$ {\it a right noetherian ring. For a
bimodule} $_{\Lambda}U _{\Gamma}$, {\it the following statements
are equivalent}.

(1) $_{\Lambda}U$ {\it is a Wakamatsu tilting module with} $\Gamma
={\rm End}(_{\Lambda}U)$.

(2) $U_{\Gamma}$ {\it is a Wakamatsu tilting module with} $\Lambda
={\rm End}(U_{\Gamma})$.

(3) $_{\Lambda}U _{\Gamma}$ {\it is a faithfully balanced
self-orthogonal bimodule}.

\vspace{0.2cm}

In the rest of this paper, we shall freely use the properties of
Wakamatsu tilting modules in Proposition 2.2 without pointing it
out explicitly.

Recall from [16] that a module $U$ in mod $\Lambda$ is called a
tilting module of projective dimension $\leq r$ if it satisfies
the following conditions:

(1) The projective dimension of $_{\Lambda}U$ is at most $r$.

(2) $_{\Lambda}U$ is self-orthogonal.

(3) There exists an exact sequence in mod $\Lambda$:
$$0 \to \Lambda \to U_0 \to U_1 \to \cdots \to U_r \to 0$$
such that each $U_i \in$add$_{\Lambda}U$ for any $0 \leq i \leq
r$.

\vspace{0.2cm}

By Proposition 2.2 and [16] Theorem 1.5, we have the following
result.

\vspace{0.2cm}

{\bf Corollary 2.3} {\it Let} $\Lambda$ {\it be a left noetherian
ring,} $\Gamma$ {\it a right noetherian ring and} $_{\Lambda}U$
{\it a Wakamatsu tilting module with} $\Gamma ={\rm
End}(_{\Lambda}U)$. {\it If the projective dimensions of}
$_{\Lambda}U$ {\it and} $U_{\Gamma}$ {\it are finite, then}
$_{\Lambda}U_{\Gamma}$ {\it is a tilting bimodule} ({\it that is,
both} $_{\Lambda}U$ {\it and} $U_{\Gamma}$ {\it are tilting}) {\it
with the left and right projective dimensions identical}.

\vspace{0.2cm}

For a module $A$ in Mod $\Lambda$ (resp. Mod $\Lambda ^{op}$), we
use {\it l.}id$_{\Lambda}(A)$ (resp. {\it r.}id$_{\Lambda}(A)$) to
denote the left (resp. right) injective dimension of $A$.

\vspace{0.2cm}

{\bf Lemma 2.4} {\it Let} $\Lambda$ {\it and} $\Gamma$ {\it be
rings and} $_{\Lambda}U_{\Gamma}$ {\it a bimodule}.

(1) {\it If} $\Gamma$ {\it is a right noetherian ring, then
r.}id$_{\Gamma}(U)=$sup$\{${\it l.}fd$_{\Gamma}({\rm
Hom}_{\Lambda}(U, E))| _{\Lambda}E$ {\it is injective}$\}$. {\it
Moreover, r.}id$_{\Gamma}(U)=${\it l.}fd$_{\Gamma}({\rm
Hom}_{\Lambda}(U, Q))$ {\it for any injective cogenerator}
$_{\Lambda}Q$ {\it for} Mod $\Lambda$.

(2) {\it If} $\Lambda$ {\it is a left noetherian ring, then
l.}id$_{\Lambda}(U)=$sup$\{${\it r.}fd$_{\Lambda}({\rm
Hom}_{\Gamma}(U, E'))|E'_{\Gamma}$ {\it is injective}$\}$. {\it
Moreover, l.}id$_{\Lambda}(U)=${\it r.}fd$_{\Lambda}({\rm
Hom}_{\Gamma}(U, Q'))$ {\it for any injective cogenerator}
$Q'_{\Gamma}$ for Mod $\Gamma ^{op}$.

\vspace{0.2cm}

{\it Proof.} (1) By [4] Chapter VI, Proposition 5.3, for any $i
\geq 1$, we have the following isomorphism:
$${\rm Tor}_{i}^{\Gamma}(B, {\rm Hom}_{\Lambda}(U, E))\cong
{\rm Hom}_{\Lambda}({\rm Ext}_{\Gamma}^{i}(B, U), E) \eqno{(1)}$$
for any $B \in$mod $\Gamma ^{op}$ and $_{\Lambda}E$ injective.

If {\it l.}fd$_{\Gamma}({\rm Hom}_{\Lambda}(U, E))\leq n(<\infty
)$ for any injective module $_{\Lambda}E$, then the isomorphism
(1) induces ${\rm Hom}_{\Lambda}({\rm Ext}_{\Gamma}^{n+1}(B, U),
E) \cong {\rm Tor}_{n+1}^{\Gamma}(B, {\rm Hom}_{\Lambda}(U,
E))=0$. Now taking $_{\Lambda}E$ as an injective cogenerator in
mod $\Lambda$, we see that Ext$_{\Gamma}^{n+1}(B, U)=0$ and {\it
r.}id$_{\Gamma}(U)\leq n$.

Conversely, if {\it r.}id$_{\Gamma}(U)=n(<\infty )$, then
Ext$_{\Gamma}^{n+1}(B, U)=0$ for any $B \in$mod $\Gamma ^{op}$ and
\linebreak Tor$_{n+1}^{\Gamma}(B, {\rm Hom}_{\Lambda}(U, E))=0$
for any injective module $_{\Lambda}E$ by the isomorphism (1).

Let $Y$ be any module in Mod $\Gamma ^{op}$. Then
$Y=\mathop{\mathrm{lim}}\limits_{\longrightarrow}Y_{\alpha}$
(where $Y_{\alpha}$ ranges over all finitely generated submodules
of $Y$). It is well known that the functor ${\rm Tor}_i$ commutes
with $\mathop{\mathrm{lim}}\limits_{\longrightarrow}$ for any $i
\geq 0$, so Tor$_{n+1}^{\Gamma}(Y, {\rm Hom}_{\Lambda}(U, E))
\cong
\mathop{\mathrm{lim}}\limits_{\longrightarrow}$Tor$_{n+1}^{\Gamma}(Y_{\alpha},
{\rm Hom}_{\Lambda}(U, E))=0$ by the above argument. This implies
that {\it l.}fd$_{\Gamma}({\rm Hom}_{\Lambda}(U, E))\leq n$.
Consequently, we conclude that the first equality holds.

The above argument in fact proves the second equality.

(2) It is similar to the proof of (1). $\blacksquare$

\vspace{0.2cm}

Let ${_{\Lambda}U _{\Gamma}}$ be a bimodule. For a module $A$ in
Mod $\Lambda$ (resp. Mod $\Gamma ^{op}$), we call
Hom$_{\Lambda}(_{\Lambda}A, {_{\Lambda}U _{\Gamma}})$ (resp.
Hom$_{\Gamma}(A_{\Gamma}, {_{\Lambda}U_{\Gamma}})$) the dual
module of $A$ with respect to $_{\Lambda}U_{\Gamma}$, and denote
either of these modules by $A^*$. For a homomorphism $f$ between
$\Lambda$-modules (resp. $\Gamma ^{op}$-modules), we put $f^*={\rm
Hom}(f, {_{\Lambda}U} _{\Gamma})$. We use $\sigma _{A}: A
\rightarrow A^{**}$ via $\sigma _{A}(x)(f)=f(x)$ for any $x \in A$
and $f \in A^*$ to denote the canonical evaluation homomorphism.
$A$ is called $U$-torsionless (resp. $U$-reflexive) if $\sigma
_{A}$ is a monomorphism (resp. an isomorphism).

\vspace{0.2cm}

{\bf Lemma 2.5} {\it Let} $\Lambda$ {\it be a left noetherian
ring}, $\Gamma$ {\it any ring and} $_{\Lambda}U_{\Gamma}$ {\it a
bimodule. If} $\Lambda ={\rm End}(U_{\Gamma})$, $U_{\Gamma}$ {\it
is self-orthogonal and r.}id$_{\Gamma}(U) \leq n$, {\it then}
$\bigoplus _{i=0}^nV_i$ {\it is an injective cogenerator for} Mod
$\Lambda$, {\it where} $V_i$ {\it is the} $(i+1)$-st {\it term in
an injective resolution of} $_{\Lambda}U$ {\it for any} $0 \leq i
\leq n$.

\vspace{0.2cm}

{\it Proof.} Let $A$ be any module in mod $\Lambda$. Since {\it
r.}id$_{\Gamma}(U) \leq n$, Ext$_{\Gamma}^i(X, U)=0$ for any $X
\in$ mod $\Gamma ^{op}$ and $i \geq n+1$. Then, by the assumption
and [13] Theorem 2.2, it is easy to see that $A$ is $U$-reflexive
provided that Ext$_{\Lambda}^i(A, U)=0$ for any $1 \leq i \leq n$.

Let $S$ be any simple $\Lambda$-module. Then Ext$_{\Lambda}^t(S,
U)\neq 0$ for some $t$ with $0 \leq t \leq n$ (Otherwise, $S$ is
$U$-reflexive by the above argument and hence $S \cong S^{**}=0$).

Let
$$0 \to {_{\Lambda}U} \to V_{0} \to
V_{1} \to \cdots \to V_{i} \to \cdots$$ be an injective resolution
of $_{\Lambda}U$. Set $W_t={\rm Im}(V_{t-1} \to V_t)$. We then get
the following exact sequences:
$${\rm Hom}_{\Lambda}(S, W_t) \to {\rm Ext}_{\Lambda}^t(S,
U) \to 0,$$
$$0 \to {\rm Hom}_{\Lambda}(S, W_t) \to {\rm Hom}_{\Lambda}(S, V_t).$$
Because Ext$_{\Lambda}^t(S, U)\neq 0$, ${\rm Hom}_{\Lambda}(S,
W_t) \neq 0$ and ${\rm Hom}_{\Lambda}(S, V_t)\neq 0$. So ${\rm
Hom}_{\Lambda}(S, \bigoplus _{i=0}^nV_i) \neq 0$ and hence
$\bigoplus _{i=0}^nV_i$ is an injective cogenerator for Mod
$\Lambda$ by [1] Proposition 18.15. $\blacksquare$

\vspace{0.2cm}

As an application to Theorem II, we have the following result.

\vspace{0.2cm}

{\bf Proposition 2.6} {\it Let} $\Lambda$ {\it and} $\Gamma$ {\it
be left and right noetherian rings and} $_{\Lambda}U$ {\it a
Wakamatsu tilting module with} $\Gamma ={\rm End}(_{\Lambda}U)$.
{\it If} $_{\Lambda}U$ {\it is} $k$-{\it Gorenstein for all} $k$
{\it and both} {\it l.}id$_{\Lambda}(U)$ {\it and
r.}id$_{\Gamma}(U)$ {\it are finite, then} {\it
l.}id$_{\Lambda}(U)=${\it r.}id$_{\Gamma}(U)$.

\vspace{0.2cm}

{\it Proof.} Assume that {\it l.}id$_{\Lambda}(U)=m<\infty$ and
{\it r.}id$_{\Gamma}(U)=n<\infty$. Since $_{\Lambda}U$ is
$k$-Gorenstein for all $k$, by Theorem II, we have that {\it
l.}fd$_{\Gamma}($Hom$_{\Lambda}(U, \bigoplus _{i=0}^mE_i)) \leq
m$, where $E_i$ is the $(i+1)$-st term in a minimal injective
resolution of $_{\Lambda}U$ for any $i \geq 0$.

By Proposition 2.2, $_{\Lambda}U _{\Gamma}$ is a faithfully
balanced self-orthogonal bimodule. If $m<n$, then, by Lemmas 2.5
and 2.4, we have that $\bigoplus _{i=0}^nE_i$ $( \cong \bigoplus
_{i=0}^mE_i$) is an injective cogenerator for Mod $\Lambda$ and
{\it l.}fd$_{\Gamma}($Hom$_{\Lambda}(U, \bigoplus
_{i=0}^mE_i))=n$, which is a contradiction. So we have that $m
\geq n$. According to the symmetry of $k$-Gorenstein modules, we
can prove $n \geq m$ similarly. $\blacksquare$

\vspace{0.2cm}

{\bf Proposition 2.7} {\it Let} $\Lambda$ {\it be a left and right
artinian ring and} $_{\Lambda}U$ {\it a Wakamatsu tilting module
with} $\Lambda ={\rm End}(_{\Lambda}U)$. {\it If} $_{\Lambda}U$
{\it is} $k$-{\it Gorenstein for all} $k$, {\it then} {\it
l.}id$_{\Lambda}(U)=${\it r.}id$_{\Lambda}(U)$.

\vspace{0.2cm}

{\it Proof.} By Theorem II, for any $i \geq 1$ and $M \in$mod
$\Lambda$ or mod $\Lambda ^{op}$, we have that
s.grade$_U$Ext$_{\Lambda}^i(M, U) \geq i$. By Proposition 2.2,
$_{\Lambda}U _{\Lambda}$ is a faithfully balanced self-orthogonal
bimodule. It then follows from [9] Theorem and its dual statement
that {\it l.}id$_{\Lambda}(U)$ is finite if and only if {\it
r.}id$_{\Lambda}(U)$ is finite. Now our conclusion follows from
Proposition 2.6. $\blacksquare$

\vspace{0.2cm}

Putting $_{\Lambda}U={_{\Lambda}\Lambda}$, we immediately have the
following result, which generalizes [2] Corollary 5.5(b).

\vspace{0.2cm}

{\bf Corollary 2.8} {\it Let} $\Lambda$ {\it be a left and right
artinian ring}. {\it If} $\Lambda$ {\it is} $k$-{\it Gorenstein
for all} $k$, {\it then} {\it l.}id$_{\Lambda}(\Lambda)=${\it
r.}id$_{\Lambda}(\Lambda)$.

\vspace{0.5cm}

\centerline{\bf 3. The proof of main results}

\vspace{0.2cm}

In this section, we prove Theorems I and II.

From now on, $\Lambda$ and $\Gamma$ are left and right noetherian
rings and $_{\Lambda}U$ is a Wakamatsu tilting module with $\Gamma
={\rm End}(_{\Lambda}U)$. We always assume that
$$0 \to {_{\Lambda}U} \to E_{0} \to
E_{1} \to \cdots \to E_{i} \to \cdots$$ is a minimal injective
resolution of $_{\Lambda}U$, and

$$0 \to U _{\Gamma} \to
E_{0}' \to E_{1}' \to \cdots \to E_{i}' \to \cdots$$ is a minimal
injective resolution of $U _{\Gamma}$ and $k$ is a positive
integer.

\vspace{0.2cm}

{\bf Lemma 3.1} {\it Let} $_{\Lambda}E$ {\it be injective. Then}
{\it l.}fd$_{\Gamma}($Hom$_{\Lambda}(U,
E))=U$-lim.dim$_{\Lambda}(E)$.

\vspace{0.2cm}

{\it Proof.} We first prove that $U$-lim.dim$_{\Lambda}(E)
\leq${\it l.}fd$_{\Gamma}($Hom$_{\Lambda}(U, E))$. Without loss of
generality, assume that {\it l.}fd$_{\Gamma}($Hom$_{\Lambda}(U,
E))=n<\infty$. Then there exists an exact sequence:

$$0 \to Q_{n} \to \cdots \to Q_{1} \to Q_{0} \to {\rm Hom}_{\Lambda}(U, E)
\to 0$$ in Mod $\Gamma$ with each $Q_{i}$ $\Gamma$-flat for any $0
\leq i \leq n$. By [4] Chapter VI, Proposition 5.3, we have that
$${\rm Tor}_{j}^{\Gamma}(U , {\rm Hom}_{\Lambda}(U, E)) \cong {\rm
Hom}_{\Lambda}({\rm Ext}_{\Gamma}^{j}(U , U), E)=0$$ for any $j
\geq 1$. Then we easily get an exact sequence:
$$0 \to U \otimes _{\Gamma}Q_{n} \to \cdots \to
U \otimes _{\Gamma}Q_{1} \to U \otimes _{\Gamma}Q_{0} \to U
\otimes _{\Gamma} {\rm Hom}_{\Lambda}(U, E) \to 0.$$ Because each
$Q_i$ is a direct limit of finitely generated free
$\Gamma$-modules, $U \otimes _{\Gamma}Q_{i}
\in$add-lim$_{\Lambda}U$ for any $0 \leq i \leq n$. On the other
hand, $U \otimes _{\Gamma} {\rm Hom}_{\Lambda}(U, E) \cong {\rm
Hom}_{\Lambda}({\rm Hom}_{\Gamma}(U, U), E) \cong E$ by [18] p.47.
So we conclude that $U$-lim.dim$_{\Lambda}(E)\leq n$.

We next prove that {\it l.}fd$_{\Gamma}($Hom$_{\Lambda}(U, E))\leq
U$-lim.dim$_{\Lambda}(E)$. Assume that
$U$-lim.dim$_{\Lambda}(E)$\linebreak $=n <\infty$. Then there
exists an exact sequence: $$0 \to X_{n} \to \cdots \to X_{1} \to
X_{0} \to E \to 0 \eqno{(2)}$$ in Mod $\Lambda$ with each $X_{i}$
in add-lim$_{\Lambda}U$ for any $0 \leq i \leq n$. Since
$_{\Lambda}U$ is finitely generated, by [17] Theorem 3.2, for any
direct system $\{{M_{\alpha}}\}_{{\alpha} \in I}$ and $j\geq 0$,
we have that Ext$_{\Lambda}^j(U,
\mathop{\mathrm{lim}}\limits_{\longrightarrow} M_{\alpha}) \cong
\mathop{\mathrm{lim}}\limits_{\longrightarrow}$Ext$_{\Lambda}^j(U,
M_{\alpha})$. From this fact we know that Ext$_{\Lambda}^{j}(U ,
X_{i})=0$ and Hom$_{\Lambda}(U , X_{i})$ is in
add-lim$_{\Gamma}\Gamma$ for any $j \geq 1$ and $0 \leq i \leq n$.
So each Hom$_{\Lambda}(U , X_{i})$ is $\Gamma$-flat for any $0
\leq i \leq n$ and by applying the functor Hom$_{\Lambda}(U, \ )$
to the exact sequence (2) we obtain the following exact sequence:
$$0 \to {\rm Hom}_{\Lambda}(U, X_n) \to \cdots \to {\rm
Hom}_{\Lambda}(U, X_1) \to {\rm Hom}_{\Lambda}(U, X_0) \to {\rm
Hom}_{\Lambda}(U, E) \to 0.$$ Hence {\it
l.}fd$_{\Gamma}($Hom$_{\Lambda}(U, E))\leq n$. The proof is
finished. $\blacksquare$

\vspace{0.2cm}

{\bf Lemma 3.2} {\it Let} $m$ {\it be an integer with} $m \geq
-k$. {\it Then the following statements are equivalent}.

(1) $U$-lim.dim$_{\Lambda}(\bigoplus _{i=0}^{k-1}E_{i})\leq k+m$.

(2) s.grade$_U{\rm Ext} _{\Gamma}^{k+m+1}(N, U) \geq k$ {\it for
any} $N \in$mod $\Gamma ^{op}$.

(3) {\it l.}fd$_{\Gamma}($Hom$_{\Lambda}(U, E_i)) \leq k+m$ {\it
for any} $0 \leq i \leq k-1$.

\vspace{0.2cm}

{\it Proof.} By Lemma 3.1, we have $(1) \Leftrightarrow (3)$.

$(2) \Rightarrow (3)$ We proceed by using induction on $i$.
Suppose that s.grade$_U{\rm Ext}_{\Gamma}^{k+m+1}(N, U)\geq k$ for
any $N \in$mod $\Gamma ^{op}$. We first prove {\it
l.}fd$_{\Gamma}($Hom$_{\Lambda}(U, E_0)) \leq k+m$. By assumption,
we have Hom$_{\Lambda}({\rm Ext}_{\Gamma}^{k+m+1}(N, U), U)=0$. We
now claim that Hom$_{\Lambda}({\rm Ext}_{\Gamma}^{k+m+1}(N, U),
E_{0})=0$. For if otherwise, then there exists $0\neq f: {\rm
Ext}_{\Gamma}^{k+m+1}(N, U)\to E_{0}$ and Im$f\bigcap U \neq 0$
(since $U$ is essential in $E_{0}$). Hence, there is a submodule
$X$($=f^{-1}({\rm Im}f\bigcap U)$) of ${\rm
Ext}_{\Gamma}^{k+m+1}(N, U)$ such that Hom$_{\Lambda}(X, U)\neq
0$, which contradicts s.grade$_U{\rm Ext}_{\Gamma}^{k+m+1}(N,
U)\geq k$. It follows easily from [4] Chapter VI, Proposition 5.3
that {\it l.}fd$_{\Gamma}({\rm Hom}_{\Lambda}(_{\Lambda}U
_{\Gamma}, E_{0})) \leq k+m$.

Now suppose $i\geq 1$. Consider the exact sequence:
$$0 \to K_{i-1}\to E_{i-1} \to K_{i} \to 0$$
where $K_{i-1}={\rm Ker}(E_{i-1}\to E_{i})$ and $K_{i}={\rm
Im}(E_{i-1}\to E_{i})$. Then for any $X \subset {\rm
Ext}_{\Gamma}^{k+m+1}(N, U)$, we have an exact sequence:
$${\rm Hom}_{\Lambda}(X, E_{i-1}) \to {\rm Hom}_{\Lambda}(X, K_{i})
\to {\rm Ext}_{\Lambda}^{1}(X, K_{i-1})\to 0 \eqno{(3)}$$ Since
s.grade$_U{\rm Ext}_{\Gamma}^{k+m+1}(N, U)\geq k$ and $1 \leq i
\leq k-1$, Ext$^{1}_{\Lambda}(X,
K_{i-1})\cong$Ext$^{i}_{\Lambda}(X, U)=0$. By induction
assumption, {\it l.}fd$_{\Gamma}($Hom$_{\Lambda}(U, E_{i-1})) \leq
k+m$. It follows from [4] Chapter VI, Proposition 5.3 that
Hom$_{\Lambda}({\rm Ext}_{\Gamma}^{k+m+1}(N, U), E_{i-1})
\cong$Tor$_{k+m+1}^{\Gamma}(N, {\rm Hom}_{\Lambda}(U,
E_{i-1}))=0$. Since $E_{i-1}$ is injective, Hom$_{\Lambda}(X,
E_{i-1})=0$. It follows from the exactness of the sequence (3)
that Hom$_{\Lambda}(X, K_{i})=0$. Observe that $E_{i}$ is the
injective envelope of $K_{i}$, by using a similar argument to the
case $i=0$, we can show that Hom$_{\Lambda}({\rm
Ext}_{\Gamma}^{k+m+1}(M, U), E_{i})=0$. Hence, we have that {\it
l.}fd$_{\Gamma}({\rm Hom}_{\Lambda}(_{\Lambda}U _{\Gamma}, E_{i}))
\leq k+m$.

$(3) \Rightarrow (2)$ Suppose that {\it l.}fd$_{\Gamma}({\rm
Hom}_{\Lambda}(U, \bigoplus _{i=0}^{k-1}E_{i}))\leq k+m$. Then, by
[4] Chapter VI, Proposition 5.3, we have that Hom$_{\Lambda}({\rm
Ext}_{\Gamma}^{k+m+1}(N, U), \bigoplus _{i=0}^{k-1}E_{i})=0$ for
any $N \in$mod $\Gamma ^{op}$. Let $X$ be any submodule of
Ext$_{\Gamma}^{k+m+1}(N, U)$. Then Hom$_{\Lambda}(X, \bigoplus
_{i=0}^{k-1}E_{i})=0$. Putting $K_{0}=U$ and $K_{i}={\rm
Im}(E_{i-1}\to E_{i})$ for any $1 \leq i \leq k-1$. Then
Hom$_{\Lambda}(X, K_{i})=0$ for any $0 \leq i \leq k-1$. It is not
difficult to prove that Ext$^{i+1}_{\Lambda}(X, K_{0})\cong$
Ext$^{1}_{\Lambda}(X, K_{i})$ and Ext$^{1}_{\Lambda}(X, K_{i})
\cong$ Hom$_{\Lambda}(X, K_{i+1})$ for any $0 \leq i \leq k-2$.
Hence we conclude that Hom$_{\Lambda}(X,
U)=0=$Ext$_{\Lambda}^{i}(X, U)$ for any $1 \leq i \leq k-1$. This
completes the proof. $\blacksquare$

\vspace{0.2cm}

Putting $m=-1$, then by Lemma 3.2, we have the following

\vspace{0.2cm}

{\bf Corollary 3.3} (1) $U$-lim.dim$_{\Lambda}(\bigoplus
_{i=0}^{k-1}E_{i})\leq k-1$ {\it if and only if} s.grade$_U{\rm
Ext} _{\Gamma}^{k}(N, U) \geq k$ {\it for any} $N\in$mod $\Gamma
^{op}$ {\it if and only if} {\it
l.}fd$_{\Gamma}($Hom$_{\Lambda}(U, \bigoplus _{i=0}^{k-1}E_{i}))
\leq k-1$.

(2) $U$-lim.dim$_{\Lambda} (E_{i})\leq i$ {\it for any} $0 \leq i
\leq k-1$ {\it if and only if} s.grade$_U{\rm Ext}
_{\Gamma}^{i}(N, U) \geq i$ {\it for any} $N\in$mod $\Gamma ^{op}$
{\it and} $1 \leq i \leq k$ {\it if and only if} {\it
l.}fd$_{\Gamma}($Hom$_{\Lambda}(U, E_i)) \leq i$ {\it for any} $0
\leq i \leq k-1$.

\vspace{0.2cm}

Let $M$ be in mod $\Lambda$ (resp. mod $\Gamma ^{op}$) and $P_{1}
\buildrel {f} \over \longrightarrow P_{0} \to M \to 0$ be a
projective presentation of $M$ in mod $\Lambda$ (resp. mod $\Gamma
^{op}$). Then we have an exact sequence:
$$0 \to M^* \to P_{0}^*
\buildrel {f^*} \over \longrightarrow P_{1}^* \to {\rm Coker}f^*
\to 0.$$ We call  ${\rm Coker}f^*$ the transpose (with respect to
$_{\Lambda}U_{\Gamma}$) of  $M$, and denote it by ${\rm Tr}_UM$.

For a positive integer $k$, recall from [10] that $M$ is called
$U$-$k$-torsionfree if Ext$_{\Gamma}^{i} ({\rm Tr}_UM, U)$ (resp.
Ext$_{\Lambda}^{i} ({\rm Tr}_UM, U))=0$ for any $1 \leq i \leq k$.
We call $M$ $U$-$k$-syzygy if there exists an exact sequence $0
\to M \to X_{0} \to X_{1} \to \cdots \buildrel {f} \over \to
X_{k-1}$ with all $X_{i}$ in add$_{\Lambda}U$ (resp.
add$U_{\Gamma}$), and denote $M$ by $\Omega ^{k}_{U}({\rm
Coker}f)$. Putting $_{\Lambda}U _{\Gamma}={_{\Lambda}\Lambda}
_{\Lambda}$, then, in this case, the notions of
$U$-$k$-torsionfree modules and $U$-$k$-syzygy modules are just
that of $k$-torsionfree modules and $k$-syzygy modules
respectively (see [3] for the definitions of $k$-torsionfree
modules and $k$-syzygy modules). We use $\cal{T}$$_{U}^{k}({\rm
mod}\ \Lambda)$ (resp. $\cal{T}$$_{U}^{k}({\rm mod}\ \Gamma
^{op})$) and $\Omega ^{k}_{U}({\rm mod}\ \Lambda)$ (resp. $\Omega
^{k}_{U}({\rm mod}\ \Gamma ^{op})$) to denote the full subcategory
of mod $\Lambda$ (resp. mod $\Gamma ^{op}$) consisting of
$U$-$k$-torsionfree modules and $U$-$k$-syzygy modules,
respcetively. It is not difficult to verify that
$\cal{T}$$_{U}^{k}({\rm mod}\ \Lambda) \subseteq \Omega
^{k}_{U}({\rm mod}\ \Lambda)$ and $\cal{T}$$_{U}^{k}({\rm mod}\
\Gamma ^{op}) \subseteq \Omega ^{k}_{U}({\rm mod}\ \Gamma ^{op})$.

The following result generalizes [3] Proposition 1.6(a).

\vspace{0.2cm}

{\bf Lemma 3.4} {\it The following statements are equivalent}.

(1) grade$_{U}$Ext$_{\Lambda}^{i+1}(M, U)\geq i$ {\it for any}
$M\in$ mod $\Lambda$ {\it and} $1\leq i \leq k-1$.

(1)$^{op}$ grade$_{U}$Ext$_{\Gamma}^{i+1}(N, U)\geq i$ {\it for
any} $N\in$ mod $\Gamma ^{op}$ {\it and} $1\leq i \leq k-1$.

{\it If one of the above equivalent conditions holds, then}
$\cal{T}$$^{i}_U({\rm mod}\ \Lambda)=\Omega ^{i}_U({\rm mod}\
\Lambda)$ {\it and} $\cal{T}$$^{i}_U({\rm mod}\ \Gamma
^{op})=\Omega ^{i}_U({\rm mod}\ \Gamma ^{op})$ {\it for any} $1
\leq i \leq k$.

\vspace{0.2cm}

{\it Proof.} The equivalence of (1) and (1)$^{op}$ was proved in
[12] Lemma 3.3. The latter assertion follows from [10] Theorem
3.1.

\vspace{0.2cm}

Putting $m=0$, then by Lemma 3.2, we have the following result, in
which the second assertion is just [3] Proposition 2.2 when
$_{\Lambda}U_{\Gamma}={_{\Lambda}\Lambda}_{\Lambda}$.

\vspace{0.2cm}

{\bf Corollary 3.5} (1) $U$-lim.dim$_{\Lambda}(\bigoplus
_{i=0}^{k-1}E_{i})\leq k$ {\it if and only if} s.grade$_U{\rm Ext}
_{\Gamma}^{k+1}(N, U) \geq k$ {\it for any} $N\in$mod $\Gamma
^{op}$ {\it if and only if l.}fd$_{\Gamma}($Hom$_{\Lambda}(U,
\bigoplus _{i=0}^{k-1}E_{i})) \leq k$.

(2) $U$-lim.dim$_{\Lambda} (E_{i})\leq i+1$ {\it for any} $0 \leq
i \leq k-1$ {\it if and only if} s.grade$_U{\rm Ext}
_{\Gamma}^{i+1}(N, U) \geq i$ {\it for any} $N\in$mod $\Gamma
^{op}$ {\it and} $1 \leq i \leq k$ {\it if and only if
l.}fd$_{\Gamma}($Hom$_{\Lambda}(U, E_i)) \leq i+1$ {\it for any}
$0 \leq i \leq k-1$. {\it In this case,} $\cal{T}$$^{i}_U({\rm
mod}\ \Gamma ^{op})=\Omega ^{i}_U({\rm mod}\ \Gamma ^{op})$ {\it
for any} $1 \leq i \leq k$.

\vspace{0.2cm}

{\it Proof.} Our assertions follows from Lemma 3.2 and Lemma 3.4.
$\blacksquare$

\vspace{0.2cm}

Putting $m=-k$, then by Lemma 3.2, we have the following

\vspace{0.2cm}

{\bf Corollary 3.6} {\it The following statements are equivalent}.

(1) $U$-dom.dim$(_{\Lambda}U) \geq k$.

(2) s.grade$_U$Ext$_{\Gamma}^1(N, U) \geq k$ {\it for any}
$N\in$mod $\Gamma ^{op}$.

(3) Hom$_{\Lambda}(U, E_i)$ {\it is} $\Gamma$-{\it flat for any}
$0 \leq i \leq k-1$.

\vspace{0.2cm}

Dually, we have the following

\vspace{0.2cm}

{\bf Corollary 3.6}$^{op}$ {\it The following statements are
equivalent}.

(1) $U$-dom.dim$(U_{\Gamma}) \geq k$.

(2) s.grade$_U$Ext$_{\Lambda}^1(M, U) \geq k$ {\it for any}
$M\in$mod $\Lambda$.

(3) Hom$_{\Gamma}(U, E'_i)$ {\it is} $\Lambda ^{op}$-{\it flat for
any} $0 \leq i \leq k-1$.

\vspace{0.2cm}

The following two results are cited from [11].

\vspace{0.2cm}

{\bf Lemma 3.7} ([11] Corollary 2.5) Hom$_{\Lambda}(U, E_0)$ {\it
is} $\Gamma$-{\it flat if and only if} Hom$_{\Gamma}(U, E'_0)$
{\it is} $\Lambda ^{op}$-{\it flat}.

\vspace{0.2cm}

{\bf Lemma 3.8} ([11] Lemma 2.6) {\it Let} $X$ be in mod $\Lambda$
({\it resp.} mod $\Gamma ^{op}$) {\it and} $n$ {\it a non-negative
integer. If} grade$_{U}X \geq n$ {\it and}
grade$_{U}$Ext$_{\Lambda}^{n}(X, U)$ ({\it resp.}
grade$_{U}$Ext$_{\Gamma}^{n}(X, U))\geq n+1$, {\it then}
grade$_{U}X \geq n+1$.

\vspace{0.2cm}

{\bf Lemma 3.9} {\it If} $U$-dom.dim$(U_{\Gamma}) \geq k$, {\it
then} $U$-dom.dim$(_{\Lambda}U) \geq k$.

\vspace{0.2cm}

{\it Proof.} When $k=1$, by Corollary 3.6$^{op}$,
Hom$_{\Lambda}(U, E'_0)$ is $\Lambda ^{op}$-flat. Then, by Lemma
3.7, Hom$_{\Lambda}(U, E_i)$ is $\Gamma$-flat. So
$U$-dom.dim$(_{\Lambda}U) \geq 1$ by Corollary 3.6.

Now suppose $k \geq 2$. By induction assumption,
$U$-dom.dim$(_{\Lambda}U) \geq k-1$. So, by Corollary 3.6, we have
that s.grade$_U$Ext$_{\Gamma}^1(N, U) \geq k-1$ for any $N\in$mod
$\Gamma ^{op}$.

Let $X$ be any submodule of Ext$_{\Gamma}^1(N, U)$. Then grade$_UX
\geq k-1$. By assumption and Corollary 3.6$^{op}$,
grade$_U$Ext$_{\Gamma}^i(X, U) \geq k$ for any $i \geq 1$. It
follows from Lemma 3.8 that grade$_UX \geq k$. So
s.grade$_U$Ext$_{\Gamma}^1(N, U) \geq k$ and hence
$U$-dom.dim$(_{\Lambda}U) \geq k$ by Corollary 3.6. $\blacksquare$

\vspace{0.2cm}

{\it Proof of Theorem I.} By Corollary 3.6 we have that $(1)
\Leftrightarrow (2)^{op} \Leftrightarrow (3)$, and by Lemma 3.9 we
have that $(1) \Rightarrow (1)^{op}$. The other implications
follow from the symmetry. $\blacksquare$

\vspace{0.2cm}

We now begin to prove Theorem II.

\vspace{0.2cm}

{\bf Lemma 3.10} ([12] Lemma 3.2) {\it If} s.grade$_{U}{\rm
Ext}^{i+1}_{\Gamma} (X, U)\geq i$ {\it for any} $X \in$mod
$\Lambda$ ({\it resp.} mod $\Gamma ^{op}$) {\it and} $1 \leq i
\leq k-1$, {\it then each} $k$-{\it syzygy module in} mod
$\Lambda$ ({\it resp.} mod $\Gamma ^{op}$) {\it is in} $\Omega
^{k}_{U}({\rm mod}\ \Lambda)$ ({\it resp.} $\Omega ^{k}_{U}({\rm
mod}\ \Gamma ^{op})$).

\vspace{0.2cm}

{\bf Theorem 3.11} {\it The following statements are equivalent}.

(1) s.grade$_U{\rm Ext} _{\Lambda}^i(M, U) \geq i$ {\it for any}
$M \in$mod $\Lambda$ {\it and} $1 \leq i \leq k$.

(2) ${\rm Ext}_{\Gamma}^i({\rm Ext}_{\Lambda}^i(\ , U), U)$ {\it
preserves monomorphisms in} mod $\Lambda$ {\it for any} $0 \leq i
\leq k-1$.

\vspace{0.2cm}

{\it Proof.} We proceed by using induction on $k$.

$(1) \Rightarrow (2)$ Let
$$0 \to X \buildrel {f} \over \to Y \to Z \to 0 \eqno{(4)}$$
be an exact sequence in mod $\Lambda$.

Suppose $k=1$. By assumption, s.grade$_U{\rm Ext}_{\Lambda}^1(Z,
U) \geq 1$. Since Coker$f^*$ is a submodule of ${\rm
Ext}_{\Lambda}^1(Z, U)$, $({\rm Coker}f^*)^*=0$ and $0 \to X^{**}
\buildrel {f^{**}} \over \to Y^{**}$ is exact.

Now suppose $k \geq 2$. From the exact sequence (4), we get an
exact sequence:
$${\rm Ext}_{\Lambda}^{k-1}(Z, U) \buildrel {\alpha} \over \to
{\rm Ext}_{\Lambda}^{k-1}(Y, U) \buildrel {\beta} \over \to {\rm
Ext}_{\Lambda}^{k-1}(X, U) \buildrel {\gamma} \over \to {\rm
Ext}_{\Lambda}^k(Z, U).$$

Set $A={\rm Im}\alpha$, $B={\rm Im}\beta$ and $C={\rm Im}\gamma$.
By (1), we have that grade$_UA \geq k-1$, grade$_UB \geq k-1$ and
grade$_UC \geq k$. Then we get the following exact sequences:
$$0 \to {\rm Ext}_{\Gamma}^{k-1}(B, U) \to
{\rm Ext}_{\Gamma}^{k-1}({\rm Ext}_{\Lambda}^{k-1}(Y, U), U),$$
$$0 \to {\rm Ext}_{\Gamma}^{k-1}({\rm Ext}_{\Lambda}^{k-1}(X, U), U)
\to {\rm Ext}_{\Gamma}^{k-1}(B, U).$$ Thus we get a composition of
monomorphisms:
$${\rm Ext}_{\Gamma}^{k-1}({\rm Ext}_{\Lambda}^{k-1}(X, U), U) \hookrightarrow
{\rm Ext}_{\Gamma}^{k-1}(B, U) \hookrightarrow {\rm
Ext}_{\Gamma}^{k-1}({\rm Ext}_{\Lambda}^{k-1}(Y, U), U),$$ which
is also a monomorphism.

$(2) \Rightarrow (1)$ Suppose $k=1$. Let $M$ be in mod $\Lambda$
and $X$ a submodule of Ext$_{\Lambda}^1(M, U)$. Because
Ext$_{\Lambda}^1(M, U)$ is in mod $\Gamma ^{op}$, $X$ is also in
mod $\Gamma ^{op}$. So there exist a positive integer $t$ and an
exact sequence:
$$0 \to U^t \buildrel {f} \over \to L \to M \to 0$$
such that the induced exact sequence:
$$L^* \buildrel {f^*} \over \to (U^t)^* \to {\rm Ext}_{\Lambda}^1(M, U)$$
has the property that $X \cong$Coker$f^*$. By assumption, $f^{**}$
is monic, $X^* \cong$Ker$f^{**}=0$. Hence we conclude that
s.grade$_U{\rm Ext}_{\Lambda}^1(M, U) \geq 1.$

Now suppose $k \geq 2$. By induction assumption, for any $M
\in$mod $\Lambda$, we have that s.grade$_U{\rm Ext}
_{\Lambda}^i(M, U) \geq i$ for any $1 \leq i \leq k-1$ and
s.grade$_U{\rm Ext} _{\Lambda}^k(M, U) \geq k-1$. By [10] Theorem
3.1, $\Omega _U^i({\rm mod}\ \Lambda)=\mathcal{T}_U^i({\rm mod}\
\Lambda)$ for any $1 \leq i \leq k$.

Let
$$\cdots \buildrel {g_{i+1}} \over \longrightarrow P_i
\buildrel {g_i} \over \to \cdots \buildrel {g_2} \over \to P_1
\buildrel {g_1} \over \to P_0 \to M \to 0$$ be a projective
resolution of $M$ in mod $\Lambda$. Notice that Coker$g_k$ is a
$(k-1)$-syzygy module in mod $\Lambda$, so it is in $\Omega
_U^{k-1}({\rm mod}\ \Lambda)$ by Lemma 3.10 and hence in
$\mathcal{T}_U^{k-1}({\rm mod}\ \Lambda)$. Thus
Ext$_{\Gamma}^i({\rm Coker}g_k^*, U)=0$ for any $1 \leq i \leq
k-1$.

Let $X$ be a submodule of ${\rm Ext} _{\Lambda}^k(M, U)$. Then
grade$_UX \geq k-1$. By [9] Lemma 2, there exists an embedding $0
\to X \to {\rm Coker}g_k^*$. By assumption, we then have an exact
sequence:
$$0 \to {\rm Ext}_{\Gamma}^{k-1}({\rm Ext}_{\Lambda}^{k-1}(X, U), U)
\to {\rm Ext}_{\Gamma}^{k-1}({\rm Ext}_{\Lambda}^{k-1}({\rm
Coker}g_k^*, U), U)=0,$$ which implies that ${\rm
Ext}_{\Gamma}^{k-1}({\rm Ext}_{\Lambda}^{k-1}(X, U), U)=0$.

On the other hand, s.grade$_U{\rm Ext}_{\Gamma}^{k-1}(X, U) \geq
k-1$ by [21] Theorem 7.5. So s.grade$_U$\linebreak${\rm
Ext}_{\Gamma}^{k-1}(X, U) \geq k$. It follows from Lemma 3.8 that
grade$_UX \geq k$ and s.grade$_U{\rm Ext}_{\Gamma}^k(M, U) \geq
k$. We are done. $\blacksquare$

\vspace{0.2cm}

{\it Proof of Theorem II.} By definition, we have $(1)
\Leftrightarrow (2)$. By Corollary 3.5(2), we have that $(3)
\Leftrightarrow (2)^{op} \Leftrightarrow (4)$. By Theorem 3.11 and
[21] Theorem 7.5, we have that $(5) \Leftrightarrow (2)
\Leftrightarrow (2)^{op}$. The other implications follow from the
symmetry. $\blacksquare$

\vspace{0.5cm}

\centerline{\bf 4. Exactness of the double dual}

\vspace{0.2cm}

As applications to the results in Section 3, we give in this
section some characterizations of $(-)^{**}$ preserving
monomorphisms and being left exact, respectively.

As an immediate consequence of Theorem II, we have the following
result, which generalizes [5] Theorem 1 and [7] Proposition 3.1.

\vspace{0.2cm}

{\bf Proposition 4.1} {\it The following statements are
equivalent}.

(1) $U$-dom.dim$(_{\Lambda}U) \geq 1$.

(2) s.grade$_U$Ext$_{\Lambda}^1(M, U) \geq 1$ {\it for any} $M
\in$mod $\Lambda$.

(3) $E_0 \in$add-lim$_{\Lambda}U$.

(4) $(\ )^{**}$ {\it preserves monomorphisms in} mod $\Lambda$.

(1)$^{op}$ $U$-dom.dim$(U_{\Gamma}) \geq 1$.

(2)$^{op}$ s.grade$_U$Ext$_{\Gamma}^1(N, U) \geq 1$ {\it for any}
$N \in$mod $\Gamma ^{op}$.

(3)$^{op}$ $E'_0 \in$add-lim$U_{\Gamma}$.

(4)$^{op}$ $(\ )^{**}$ {\it preserves monomorphisms in} mod
$\Gamma ^{op}$.

\vspace{0.2cm}

{\bf Lemma 4.2} {\it Assume that} $U$-dom.dim$(_{\Lambda}U)\geq
k$. {\it Then, for a module} $M$ {\it in} mod $\Lambda$,
grade$_{U}M \geq k$ {\it if} $M^*=0$.

\vspace{0.2cm}

{\it Proof.} For any $M\in$mod $\Lambda$ and $i \geq 1$, we have
an exact sequence
$${\rm Hom}_{\Lambda}(M, E_{i-1})\to {\rm Hom}_{\Lambda}(M, K_i)
\to {\rm Ext}_{\Lambda}^{i}(M, U) \to 0 \eqno{(5)}$$ where
$K_i={\rm Im}(E_{i-1} \to E_i)$.

Suppose $U$-dom.dim$(_{\Lambda}U )\geq k$. Then each $E_{i}$ is in
add-lim$_{\Lambda}U$ for any $0 \leq i \leq k-1$. So, for a given
$M \in$mod $\Lambda$ with $M^*=0$, we have that Hom$_{\Lambda}(M,
E_{i})=0$ by [17] Theorem 3.2 and Hom$_{\Lambda}(M, K_i)=0$ for
any $0 \leq i \leq k-1$. Then by the exactness of the sequence
(5), Ext$_{\Lambda}^{i}(M, U)=0$ for any $1 \leq i \leq k-1$, and
so grade$_{U}M \geq k$. $\blacksquare$

\vspace{0.2cm}

{\bf Lemma 4.3} {\it If} $[{\rm Ext}_{\Lambda}^1(M, U)]^*=0$ {\it
for any} $M \in$mod $\Lambda$, {\it then} $N^*$ {\it is} $U$-{\it
reflexive for any} $N \in$mod $\Gamma ^{op}$.

\vspace{0.2cm}

{\it Proof.} By the dual statements of [10] Proposition 4.2 and
Corollary 4.2. $\blacksquare$

\vspace{0.2cm}

We now characterize $U$-dominant dimension of $U$ at least two.
The following result generalizes [5] Theorem 2 and [8] Proposition
E.

\vspace{0.2cm}

{\bf Proposition 4.4} {\it The following statements are
equivalent.}

(1) $U$-dom.dim$(_{\Lambda}U) \geq 2$.

(2) $(-)^{**}:$ mod $\Lambda \to$ mod $\Lambda$ {\it is left
exact}.

(3) $(-)^{**}:$ mod $\Lambda \to$ mod $\Lambda$ {\it preserves
monomorphisms and} Ext$_{\Gamma}^{1}({\rm Ext} _{\Lambda}^{1}(X,
U), U)=0$ {\it for any} $X\in$mod $\Lambda$.

(1)$^{op}$ $U$-dom.dim$(U _{\Gamma}) \geq 2$.

(2)$^{op}$ $(-)^{**}:$ mod $\Gamma ^{op} \to$ mod $\Gamma ^{op}$
{\it is left exact}.

(3)$^{op}$ $(-)^{**}:$ mod $\Gamma ^{op}\to$ mod $\Gamma ^{op}$
{\it preserves monomorphisms and} Ext$_{\Lambda}^{1}({\rm Ext}
_{\Gamma}^{1}(Y, U), U)=0$ {\it for any} $Y\in$mod $\Gamma ^{op}$.

\vspace{0.2cm}

{\it Proof.} By Theorem I, we have $(1)\Leftrightarrow (1)^{op}$.
By symmetry, we only need to prove that $(1) \Rightarrow (2)$ and
$(2)^{op} \Rightarrow (3) \Rightarrow (1)^{op}$.

$(1) \Rightarrow (2)$ Assume that $U$-dom.dim$(_{\Lambda}U) \geq
2$ and $0 \to A \buildrel {\alpha} \over \longrightarrow B
\buildrel {\beta} \over \longrightarrow C \to 0$ is an exact
sequence in mod $\Lambda$. By Proposition 4.1, $\alpha ^{**}$ is
monic. By Theorem I and [4] Chapter VI, Proposition 5.3, we have
that Hom$_{\Gamma}(U, E'_{0})$ is $\Lambda ^{op}$-flat and
Hom$_{\Gamma}({\rm Ext}_{\Lambda}^{1}(C, U), E'_{0})=0$. Since
Coker$\alpha ^*$ is isomorphic to a submodule of
Ext$_{\Lambda}^{1}(C, U)$, Hom$_{\Gamma} ({\rm Coker}\alpha ^*,
E'_{0})=0$ and $({\rm Coker}\alpha ^*)^*=0$. Then by Lemma 4.2, we
have that grade$_U{\rm Coker}\alpha ^* \geq 2$ and
Ext$^1_{\Gamma}({\rm Coker}\alpha ^*, U)=0$. It follows easily
that $0 \to A^{**} \buildrel {\alpha ^{**}} \over \longrightarrow
B^{**} \buildrel {\beta ^{**}} \over \longrightarrow C^{**}$ is
exact.

$(2)^{op} \Rightarrow (3)$ By Proposition 4.1, $(-)^{**}:$ mod
$\Gamma ^{op} \to$ mod $\Gamma ^{op}$ preserves monomorphisms and
$U$-dom.dim$(_{\Lambda}U)=U$-dom.dim$(U_{\Gamma}) \geq 1$. By
Theorem I, for any $X \in$ mod $\Lambda$, we have that
s.grade$_U{\rm Ext}_{\Lambda}^1(X, U) \geq 1$ and $[{\rm
Ext}_{\Lambda}^1(X, U)]^*=0$.

Let
$$0 \to K \buildrel {f} \over \to Q \buildrel {g} \over \to
{\rm Ext}_{\Lambda}^1(X, U) \to 0$$ be an exact sequence in mod
$\Gamma ^{op}$ with $Q$ projective. Then, by $(2)^{op}$, $f^{**}$
and $f^{***}$ are monomorphisms and hence isomorphisms. On the
other hand, we have the following commutative diagram with exact
rows:
$$\xymatrix{0 \ar[r] & Q^* \ar[r]^{f^*} \ar[d]^{\sigma _{Q^*}} &
K^* \ar[r] \ar[d]^{\sigma _{K^*}} & {\rm Ext}_{\Gamma}^1({\rm
Ext}_{\Lambda}^1(X, U), U) \ar[r] & 0 &\\
& Q^{***} \ar[r]^{f^{***}} & K^{***} & & &\\ }$$ It follows from
Lemma 4.3 that $Q^*$ and $K^*$ are $U$-reflexive. So $\sigma
_{Q^*}$ and $\sigma _{K^*}$ are isomorphisms and hence $f^*$ is an
isomorphism. Consequently we have that ${\rm Ext}_{\Gamma}^1({\rm
Ext}_{\Lambda}^1(X, U), U)=0$.

$(3)\Rightarrow (1)^{op}$ Suppose that (3) holds. Then
$U$-dom.dim$(U _{\Gamma}) \geq 1$ by Proposition 4.1.

Let $A$ be in mod $\Lambda$ and $B$ any submodule of
Ext$_{\Lambda}^{1}(A, U)$ in mod $\Gamma ^{op}$. Since
$U$-dom.dim$(U _{\Gamma})$\linebreak $\geq 1$, by Theorem I and
[4] Chapter VI, Proposition 5.3, we have that Hom$_{\Gamma}(U,
E'_{0})$ is $\Lambda ^{op}$-flat and Hom$_{\Gamma}({\rm
Ext}_{\Lambda}^{1}(A, U), E_{0}^{'})=0$. So Hom$_{\Gamma}(B,
E_{0}^{'})=0$ and hence Hom$_{\Gamma}(B, E_{0}^{'}/U
_{\Gamma})\cong {\rm Ext}_{\Gamma}^{1}(B, U _{\Gamma})$. On the
other hand, Hom$_{\Gamma}(B, E_{0}^{'})=0$ implies $B^*=0$. Then
by [13] Lemma 2.1, we have that $B\cong$Ext$_{\Lambda}^{1}({\rm
Tr}_{U}B, U)$ with ${\rm Tr}_{U}B$ in mod $\Lambda$. By (3),
Hom$_{\Gamma}(B, E_{0}^{'}/U)\cong {\rm Ext}_{\Gamma}^{1}(B, U)
\cong {\rm Ext}_{\Gamma}^{1}({\rm Ext}_{\Lambda}^{1}({\rm
Tr}_{U}B, U), U)=0$. Then by using a similar argument to the proof
of $(2) \Rightarrow (3)$ in Lemma 3.2, we have that
Hom$_{\Gamma}({\rm Ext}_{\Lambda}^{1}(A, U), E_{1}^{'})=0$ (note:
$E_{1}^{'}$ is the injective envelope of $E_{0}^{'}/U$). It
follows from [4] Chapter VI, proposition 5.3 that Hom$_{\Gamma}(U,
E_{1}^{'})$ is $\Lambda ^{op}$-flat and thus $U$-dom.dim$(U
_{\Gamma}) \geq 2$ by Theorem I. $\blacksquare$

\vspace{0.5cm}

\centerline{\bf 5. A generalization of $k$-Gorenstein modules}

\vspace{0.2cm}

In this section, we study a generalization of $k$-Gorenstein
modules, which is however not left-right symmetric. We
characterize this generalization in terms of some properties
similar to that of $k$-Gorenstein modules. The results obtained
here develops the main result of Auslander and Reiten in [3].

We begin with the following equivalent characterizations of
$U$-lim.dim$_{\Lambda}(E_{0}) \leq 1$ as follows, which
generalizes [8] Proposition D.

\vspace{0.2cm}

{\bf Proposition 5.1} {\it The following statements are
equivalent.}

(1) $U$-lim.dim$_{\Lambda}(E_{0}) \leq 1$.

(2) $\sigma _{X}$ {\it is an essential monomorphism for any}
$U$-{\it torsionless module} $X$ {\it in} mod $\Lambda$.

(3) $f^{**}$ {\it is a monomorphism for any monomorphism} $f:X\to
Y$ {\it in} mod $\Lambda$ {\it with} $Y$ $U$-{\it torsionless}.

(4) $f^{**}$ {\it is a monomorphism for any monomorphism} $f:X\to
Y$ {\it in} mod $\Lambda$ {\it with} $X$ {\it and} $Y$ $U$-{\it
torsionless}.

(5) grade$_{U}$Ext$^1_{\Lambda}(X, U) \geq 1$ {\it for any} $X
\in$mod $\Lambda$.

(6) s.grade$_{U}$Ext$^2_{\Gamma}(N, U) \geq 1$ {\it for any}
$N\in$mod $\Gamma ^{op}$.

\vspace{0.2cm}

{\it Proof.} $(1) \Leftrightarrow (6)$ follows from Corollary
3.5(2) and $(3) \Rightarrow (4)$ is trivial.

$(1) \Rightarrow (2)$ Suppose $U$-lim.dim$_{\Lambda}(E_{0}) \leq
1$. Then by Lemma 3.1, we have that {\it
l.}fd$_{\Gamma}($Hom$_{\Lambda}(U,$\linebreak $E_0)) \leq 1$.

Assume that $X$ is $U$-torsionless in mod $\Lambda$. Then
Coker$\sigma _{X}\cong {\rm Ext}_{\Gamma}^{2}({\rm Tr}_{U}X, U)$
by [13] Lemma 2.1. By [4] Chapter VI, Proposition 5.3, we have
that Hom$_{\Lambda}($Coker$\sigma _{X}, E_{0})\cong$\linebreak
Hom$_{\Lambda}($Ext$_{\Gamma}^{2}($Tr$_{U}X, U),
E_{0})\cong$Tor$_2^{\Gamma}({\rm Tr}_{U}X, {\rm Hom}_{\Lambda}(U,
E_0))=0$. Then $A^*=0$ for any submodule $A$ of Coker$\sigma
_{X}$, which implies that any non-zero submodule of Coker$\sigma
_{X}$ is not $U$-torsionless.

Let $B$ be a submodule of $X^{**}$ with $X \bigcap B=0$. Then $B
\cong B/(X\bigcap B) \cong (X+B)/X$ is isomorphic to a submodule
of Coker$\sigma _{X}$. On the other hand, $B$ is clearly
$U$-torsionless. So $B=0$ and hence $\sigma _{X}$ is essential.

$(2) \Rightarrow (3)$ Let $f: X\to Y$ be monic in mod $\Lambda$
with $Y$ $U$-torsionless. Then $f^{**}\sigma _{X}=\sigma _{Y}f$ is
monic. By (2), $\sigma _{X}$ is an essential monomorphism, so
$f^{**}$ is monic.

$(4) \Rightarrow (5)$ Let $X$ be in mod $\Lambda$ and $0 \to Y
\buildrel {g} \over \longrightarrow P \to X \to 0$ an exact
sequence in mod $\Lambda$ with $P$ projective. It is easy to see
that $[{\rm Ext}^{1}_{\Lambda}(X, U)]^* \cong {\rm Ker}g^{**}$. On
the other hand, since $_{\Lambda}U_{\Gamma}$ is a faithfully
balanced bimodule, $P$ is $U$-reflexive and $Y$ is
$U$-torsionless. So $g^{**}$ is monic by (4) and hence
Ker$g^{**}=0$ and $[{\rm Ext}^{1}_{\Lambda}(X, U)]^*=0$.

$(5) \Rightarrow (1)$ Let $M$ be in mod $\Gamma ^{op}$ and $\cdots
\to P_{1} \to P_{0} \to M \to 0$ a projective resolution of $M$ in
mod $\Gamma ^{op}$. Put $N=$Coker$(P_{2}\to P_{1})$. By [13] Lemma
2.1, ${\rm Ext}^{2}_{\Gamma}(M, U) \cong {\rm Ext}^{1}_{\Gamma}(N,
U) \cong {\rm Ker}\sigma _{{\rm Tr}_{U}N}$. On the other hand,
since $N$ is $U$-torsionless, ${\rm Ext}^{1}_{\Lambda}({\rm
Tr}_{U}N, U)\cong {\rm Ker}\sigma _{N}=0$.

Let $X$ be any finitely generated submodule of ${\rm
Ext}^{2}_{\Gamma}(M, U)$ and $f_{1}: X\to{\rm Ext}^{2}_{\Gamma}(M,
U) (\cong {\rm Ker}\sigma _{{\rm Tr}_{U}N})$ the inclusion, and
let $f$ be the composition: $X \buildrel {f_{1}} \over
\longrightarrow {\rm Ext}^{2}_{\Gamma}(M, U) \buildrel {g} \over
\longrightarrow {\rm Tr}_{U}N$, where $g$ is a monomorphism. Then
$\sigma _{{\rm Tr}_{U}N}f=0$ and $f^*\sigma _{{\rm
Tr}_{U}N}^*=(\sigma _{{\rm Tr}_{U}N}f)^*=0$. But $\sigma _{{\rm
Tr}_{U}N}^*$ is epic by [1] Proposition 20.14, so $f^*=0$. Hence,
by applying the functor Hom$_{\Lambda}(-, U)$ to the exact
sequence $0 \to X \buildrel {f} \over \longrightarrow {\rm
Tr}_{U}N \to {\rm Coker}f \to 0$, we have that $X^* \cong {\rm
Ext}^{1}_{\Lambda}({\rm Coker}f, U)$ and then $X^{**} \cong [{\rm
Ext}^{1}_{\Lambda}({\rm Coker}f, U)]^*=0$ by (5), which implies
that $X^*=0$ since $X^*$ is a direct summand of $X^{***}(=0)$. By
using a similar argument to the proof of $(2) \Rightarrow (3)$ in
Lemma 3.2, we can prove that {\it
l.}fd$_{\Gamma}($Hom$_{\Lambda}(U, E_0)) \leq 1$. Therefore
$U$-lim.dim$_{\Lambda}(E_{0}) \leq 1$ by Lemma 3.1. $\blacksquare$

\vspace{0.2cm}

By Proposition 4.1, we have that $E_0 \in$add-lim$_{\Lambda}U$ if
and only if $E'_0 \in$add-lim$U_{\Gamma}$, that is,
$U$-lim.dim$_{\Lambda}(E_{0})=0$ if and only if
$U$-lim.dim$_{\Gamma}(E'_{0})=0$. However, in general, we don't
have the fact that $U$-lim.dim$_{\Lambda}(E_{0}) \leq 1$ if and
only if $U$-lim.dim$_{\Gamma}(E'_{0}) \leq 1$ even when
$_{\Lambda}U_{\Gamma}={_{\Lambda}\Lambda _{\Lambda}}$.

\vspace{0.2cm}

{\bf Example} We use $I_0$ and $I'_0$ to denote the injective
envelope of $_{\Lambda}\Lambda$ and $\Lambda _{\Lambda}$,
respectively. Consider the following example. Let $K$ be a field
and $\Delta$ the quiver:
$$\xymatrix{ 1 \ar @<2pt> [r]^{\alpha} &
2 \ar @<2pt> [l]^{\beta} \ar[r]^{\gamma} & 3 }
$$
(1) If $\Lambda =K\Delta /(\alpha \beta \alpha)$. Then {\it
l.}fd$_{\Lambda}(I_0) =1$ and {\it r.}fd$_{\Lambda}(I'_0) \geq 2$.
(2) If $\Lambda =K\Delta /(\gamma \alpha , \beta \alpha)$. Then
{\it l.}fd$_{\Lambda}(I_0) =2$ and {\it r.}fd$_{\Lambda}(I'_0)
=1$.

\vspace{0.2cm}

Compare the following result with Theorem 3.11.

\vspace{0.2cm}

{\bf Theorem 5.2} {\it The following statements are equivalent}.

(1) grade$_U{\rm Ext} _{\Lambda}^i(M, U) \geq i$ {\it for any} $M
\in$mod $\Lambda$ {\it and} $1 \leq i \leq k$.

(2) ${\rm Ext}_{\Gamma}^i({\rm Ext}_{\Lambda}^i(\ , U), U)$ {\it
preserves monomorphisms} $X \to Y$ {\it with both} $X$ {\it and} $Y$
$U$-{\it torsionless in} mod $\Lambda$ {\it for any} $0 \leq i \leq
k-1$.

\vspace{0.2cm}

{\it Proof.} We proceed by using induction on $k$. The case $k=1$
follows from Proposition 5.1. Now suppose $k \geq 2$.

$(1) \Rightarrow (2)$ Let $A$ be a $U$-torsionless module in mod
$\Lambda$. Then there exists an exact sequence in mod $\Lambda$ with
$P$ in add$_{\Lambda}U$:
$$0 \to A \to P \to B \to 0.$$ By (1),
for any $1 \leq i \leq k-1$, we have that grade$_U{\rm Ext}
_{\Lambda}^i(A, U)=$grade$_U{\rm Ext} _{\Lambda}^{i+1}(B, U) \geq
i+1$, which implies that ${\rm Ext}_{\Gamma}^i({\rm
Ext}_{\Lambda}^i(A, U), U)=0$. The desired conclusion follows
trivially.

$(2) \Rightarrow (1)$ By induction assumption, for any $M \in$mod
$\Lambda$, we have that grade$_U{\rm Ext} _{\Lambda}^i(M,
U)$\linebreak $\geq i$ for any $1 \leq i \leq k-1$ and
grade$_U{\rm Ext} _{\Lambda}^k(M, U) \geq k-1$. So it suffices to
prove that ${\rm Ext}_{\Gamma}^{k-1}({\rm Ext} _{\Lambda}^k(M, U),
U)=0$.

Let
$$0 \to K \to P \to M \to 0$$ be an exact sequence in mod
$\Lambda$ with $P$ projective. Then by (2), we have the following
exact sequence:
$$0 \to {\rm Ext}_{\Gamma}^{k-1}({\rm Ext}_{\Lambda}^{k-1}(K, U),
U) \to {\rm Ext}_{\Gamma}^{k-1}({\rm Ext}_{\Lambda}^{k-1}(P, U),
U).$$ But the last term in this sequence is always zero, so ${\rm
Ext}_{\Gamma}^{k-1}({\rm Ext}_{\Lambda}^k(M, U), U) \cong {\rm
Ext}_{\Gamma}^{k-1}({\rm Ext}_{\Lambda}^{k-1}$\linebreak$(K, U),
U)=0$. $\blacksquare$

\vspace{0.2cm}

Compare the following result with [21] Theorem 7.5.

\vspace{0.2cm}

{\bf Theorem 5.3} {\it The following statements are equivalent}.

(1) s.grade$_U{\rm Ext} _{\Gamma}^{i+1}(N, U) \geq i$ {\it for
any} $N \in$mod $\Gamma ^{op}$ {\it and} $1 \leq i \leq k$.

(2) grade$_U{\rm Ext} _{\Lambda}^i(M, U) \geq i$ {\it for any} $M
\in$mod $\Lambda$ {\it and} $1 \leq i \leq k$.

\vspace{0.2cm}

{\it Proof.} We proceed by using induction on $k$. The case $k=1$
follows from Proposition 5.1. Now suppose $k \geq 2$.

$(1) \Rightarrow (2)$ By induction assumption, for any $M \in$mod
$\Lambda$, we have that grade$_U{\rm Ext} _{\Lambda}^i(M,
U)$\linebreak $\geq i$ for any $1 \leq i \leq k-1$ and
grade$_U{\rm Ext} _{\Lambda}^k(M, U) \geq k-1$. Then
$\cal{T}$$^{i}_U({\rm mod}\ \Lambda)=\Omega ^{i}_U({\rm mod}\
\Lambda)$ for any $1 \leq i \leq k$ by Lemma 3.4.

Let $$\cdots \to P_i \to \cdots \to P_1 \to P_0 \to M \to 0$$ be
an exact sequence in mod $\Lambda$ with each $P_i$ projective for
any $i \geq 0$. By [9] Lemma 2, we have the following exact
sequence:
$$0 \to {\rm Ext}_{\Lambda}^k(M, U) \to {\rm Tr}_U\Omega _{\Lambda}^{k-1}(M)
\to P_{k+1}^* \to {\rm Tr}_U\Omega _{\Lambda}^k(M) \to 0
\eqno{(6)}$$ Notice that $\Omega _{\Lambda}^{k-1}(M)$ is
$(k-1)$-syzygy and $\Omega _{\Lambda}^k(M)$ is $k$-syzygy, so, by
Lemma 3.10, $\Omega _{\Lambda}^{k-1}(M)$ (resp. $\Omega
_{\Lambda}^k(M)$) is in $\Omega _U^{k-1}({\rm mod} \ \Lambda)$
(resp. $\Omega _U^k({\rm mod} \ \Lambda)$) and hence is in
$\mathcal{T}_U^{k-1}({\rm mod} \ \Lambda)$ (resp.
$\mathcal{T}_U^k({\rm mod} \ \Lambda)$). It follows that
Ext$_{\Gamma}^i({\rm Tr}_U\Omega _{\Lambda}^{k-1}(M), U)=0$ for
any $1 \leq i \leq k-1$ and Ext$_{\Gamma}^i({\rm Tr}_U\Omega
_{\Lambda}^k(M), U)=0$ for any $1 \leq i \leq k$. In addition,
$P_{k+1}^* \in$add$U_{\Gamma}$, so Ext$_{\Gamma}^i(P_{k+1}^*,
U)=0$ for any $i \geq 1$. Thus from the exact sequence (6) we get
an embedding:
$$0 \to {\rm Ext}_{\Gamma}^{k-1}({\rm
Ext}_{\Lambda}^k(M, U), U) \to {\rm Ext}_{\Gamma}^{k+1}({\rm
Tr}_U\Omega _{\Lambda}^k(M), U).$$ Then, by (1), we have that
grade$_U{\rm Ext}_{\Gamma}^{k-1}({\rm Ext}_{\Lambda}^k(M, U), U)
\geq k$. Consequently, grade$_U{\rm Ext}_{\Lambda}^k$\linebreak
$(M, U) \geq k$ by Lemma 3.8.

$(2) \Rightarrow (1)$ By induction assumption, for any $N \in$mod
$\Gamma ^{op}$, we have that s.grade$_U{\rm Ext}
_{\Gamma}^{i+1}(N,$\linebreak $U) \geq i$ for any $1 \leq i \leq
k-1$ and s.grade$_U{\rm Ext} _{\Gamma}^{k+1}(N, U) \geq k-1$. Then
$\cal{T}$$^{i}_U({\rm mod}\ \Gamma ^{op})=\Omega ^{i}_U({\rm mod}\
\Gamma ^{op})$ for any $1 \leq i \leq k$ by Lemma 3.4.

Let $X$ be a submodule of ${\rm Ext} _{\Gamma}^{k+1}(N, U)$. Then
grade$_UX \geq k-1$. By [9] Lemma 2, there exists an exact
sequence:
$$0 \to X \buildrel {f} \over \to {\rm Tr}_U\Omega _{\Gamma}^k(N)
\to {\rm Coker}f \to 0 \eqno{(7)}$$ Notice that $\Omega
_{\Gamma}^k(N)$ is $k$-syzygy, so, by Lemma 3.10, it is in $\Omega
_U^k({\rm mod} \ \Gamma ^{op})$ and hence is in
$\mathcal{T}_U^k({\rm mod} \ \Gamma ^{op})$. It follows that
Ext$_{\Lambda}^i({\rm Tr}_U\Omega _{\Gamma}^k(N), U)=0$ for any $1
\leq i \leq k$. So from the exact sequence (7) we get that
Ext$_{\Lambda}^{k-1}(X, U) \cong$Ext$_{\Lambda}^k({\rm Coker}f,
U)$. By (2), grade$_U{\rm Ext} _{\Lambda}^{k-1}(X, U)=$\linebreak
grade$_U{\rm Ext} _{\Lambda}^k({\rm Coker}f, U) \geq k$. It
follows from Lemma 3.8 that grade$_UX \geq k$ and
s.grade$_U$\linebreak ${\rm Ext} _{\Gamma}^{k+1}(N, U) \geq k$.
$\blacksquare$

\vspace{0.2cm}

Recall that a full subcategory $\mathcal{X}$ of mod $\Lambda$
(resp. mod $\Gamma ^{op}$) is said to be closed under extensions
if the middle term $B$ of any short sequence $0 \to A \to B \to C
\to 0$ is in $\mathcal{X}$ provided that the end terms $A$ and $C$
are in $\mathcal{X}$.

The following is the main result in this section.

\vspace{0.2cm}

{\bf Theorem 5.4} {\it The following statements are equivalent}.

(1) s.grade$_U$Ext$_{\Gamma}^{i+1}(N, U) \geq i$ {\it for any} $N
\in$mod $\Gamma ^{op}$ {\it and} $1 \leq i \leq k$.

(2) $U$-lim.dim$_{\Lambda}(E_i) \leq i+1$ {\it for any} $0 \leq i
\leq k-1$.

(3) {\it l.}fd$_{\Gamma}($Hom$_{\Lambda}(U, E_i)) \leq i+1$ {\it
for any} $0 \leq i \leq k-1$ {\it for any} $0 \leq i \leq k-1$.

(4) grade$_U$Ext$_{\Lambda}^i(M, U) \geq i$ {\it for any} $M
\in$mod $\Lambda$ {\it and} $1 \leq i \leq k$.

(5) ${\rm Ext}_{\Gamma}^i({\rm Ext}_{\Lambda}^i(\ , U), U)$ {\it
preserves monomorphisms} $X \to Y$ {\it with both} $X$ {\it and} $Y$
$U$-{\it torsionless in} mod $\Lambda$ {\it for any} $0 \leq i \leq
k-1$.

{\it If one of the above equivalent conditions holds, then}
$\Omega _U^i({\rm mod} \ \Gamma ^{op})(=\mathcal{T}_U^i({\rm mod}
\ \Gamma ^{op}))$ {\it is closed under extensions for any} $1 \leq
i \leq k$.

\vspace{0.2cm}

{\it Proof.} By Corollary 3.5(2), we have that $(1)\Leftrightarrow
(2)\Leftrightarrow (3)$. It follows from Theorems 5.3 and 5.2 that
$(1)\Leftrightarrow (4)\Leftrightarrow (5)$. The last assertion
follows from [10] Theorem 3.3. $\blacksquare$

\vspace{0.2cm}

We use $I_i$ (resp. $I'_i$) to denote the $(i+1)$-st term in a
minimal injective resolution of $_{\Lambda}\Lambda$ (resp.
$\Lambda _{\Lambda}$) for any $i \geq 0$. The following corollary
generalizes [3] Theorem 4.7. In [3], the assumption of $\Lambda$
being a noetherian algebra is necessary for proving
$(5)\Rightarrow (3)$. But here the assumption of $\Lambda$ being a
left and right noetherian ring is enough for all of the
implications.

\vspace{0.2cm}

{\bf Corollary 5.5} {\it The following statements are equivalent}.

(1) s.grade$_\Lambda$Ext$_{\Lambda}^{i+1}(N, \Lambda) \geq i$ {\it
for any} $N \in$mod $\Lambda ^{op}$ {\it and} $1 \leq i \leq k$.

(2) {\it l.}fd$_{\Lambda}(I_i) \leq i+1$ {\it for any} $0 \leq i
\leq k-1$.

(3) grade$_\Lambda$Ext$_{\Lambda}^i(M, \Lambda) \geq i$ {\it for
any} $M \in$mod $\Lambda$ {\it and} $1 \leq i \leq k$.

(4) ${\rm Ext}_{\Lambda}^i({\rm Ext}_{\Lambda}^i(\ , \Lambda),
\Lambda)$ {\it preserves such monomorphisms} $X \to Y$ {\it with
both} $X$ {\it and} $Y$ {\it torsionless in} mod $\Lambda$ {\it
for any} $0 \leq i \leq k-1$.

(5) $\Omega _{\Lambda}^i({\rm mod} \ \Lambda ^{op})$ {\it is
closed under extensions for any} $1 \leq i \leq k$.

(6) add$\Omega _{\Lambda}^i({\rm mod} \ \Lambda ^{op})$ ({\it the
subcategory of} mod $\Lambda ^{op}$ {\it whose objects are those
modules which are direct summands of} $i$-th {\it syzygies}) {\it
is closed under extensions for any} $1 \leq i \leq k$.

\vspace{0.2cm}

{\it Proof.} By Theorem 5.4, we have that $(1)\Leftrightarrow
(2)\Leftrightarrow (3)\Leftrightarrow (4)$. The equivalence of
(1), (5) and (6) follows from the dual statements of [3] Theorem
1.7. $\blacksquare$

\vspace{0.2cm}

At the end of this section, we generalize the result of Wakamatsu
on the symmetry of $k$-Gorenstein modules.

\vspace{0.2cm}

{\bf Proposition 5.6} {\it Assume that} $m$ {\it is a non-negative
integer and} $U$-lim.dim$_{\Lambda}(E_i) \leq i+1$ {\it for any}
$0 \leq i \leq m-1$.

(1) {\it If} $U$-lim.dim$_{\Gamma}(\bigoplus _{i=0}^mE'_i) \leq
m$, {\it then} $U$-lim.dim$_{\Lambda}(E_m) \leq m$; {\it
Especially, if} {\it l.}id$_{\Lambda}(U) \leq m$, {\it then}
$U$-lim.dim$_{\Lambda}(E_m) \leq m$.

(2) {\it For a positive integer} $k$, {\it if}
$U$-lim.dim$_{\Gamma}(\bigoplus _{i=0}^mE'_i) \leq m$ {\it and}
$U$-lim.dim$_{\Gamma}(E'_{m+j}) \leq m+j$ {\it for any} $1 \leq j
\leq k-1$, {\it then} $U$-lim.dim$_{\Lambda}(E_{m+j}) \leq m+j$
{\it for any} $0 \leq j \leq k-1$.

\vspace{0.2cm}

{\it Proof.} The case $m=0$ follows from Theorem II. Now suppose
$m \geq 1$.

(1) By Corollaries 3.5 and 3.3, it suffices to prove that if
s.grade$_U$Ext$_{\Gamma}^{i+1}(N, U) \geq i$ for any $N \in$mod
$\Gamma ^{op}$ and $1 \leq i \leq m$ and
s.grade$_U$Ext$_{\Lambda}^{m+1}(M, U) \geq m+1$ for any $M \in$mod
$\Lambda$, then s.grade$_U$Ext$_{\Gamma}^{m+1}(N, U) \geq m+1$ for
any $N \in$mod $\Gamma ^{op}$.

Suppose that
$$\cdots \to Q_i \to \cdots \to Q_1 \to Q_0 \to N \to 0 \eqno{(8)}$$ is a
projective resolution of $N$ in mod $\Gamma ^{op}$.

By Lemma 3.4, we have that $\mathcal{T}_U^i({\rm mod} \ \Gamma
^{op})=\Omega _U^i({\rm mod} \ \Gamma ^{op})$ for any $1 \leq i
\leq m+1$. Notice that Coker$(Q_{m+1} \to Q_m)$ is $m$-syzygy, so,
by Lemma 3.10, it is in $\Omega _U^m({\rm mod} \ \Gamma ^{op})$
and hence is in $\mathcal{T}_U^m({\rm mod} \ \Gamma ^{op})$, which
implies that Ext$_{\Lambda}^i({\rm Tr}_U\Omega _{\Gamma}^m(N),
U)=0$ for any $1 \leq i \leq m$.

Let $X$ be a submodule of Ext$_{\Gamma}^{m+1}(N, U)$. Then
grade$_UX \geq m$. By [9] Lemma 2, we have an exact sequence:
$$0 \to X \buildrel {f} \over \to {\rm Tr}_U\Omega
_{\Gamma}^m(N) \to {\rm Coker}f \to 0.$$ We then get an embedding
$0 \to {\rm Ext}_{\Lambda}^m(X, U) \to {\rm
Ext}_{\Lambda}^{m+1}({\rm Coker}f, U)$. By assumption,
s.grade$_U$Ext$_{\Lambda}^{m+1}({\rm Coker}f, U) \geq m+1$. So
grade$_U$Ext$_{\Lambda}^m(X, U) \geq m+1$. It follows from Lemma
3.8 that grade$_UX \geq m+1$ and s.grade$_U$Ext$_{\Gamma}^{m+1}(N,
U) \geq m+1$.

By Lemma 2.4(2) and the dual statement of Lemma 3.1, we have that
$U$-lim.dim$_{\Gamma}$\linebreak $(\bigoplus _{i=0}^kE'_i)
\leq${\it l.}id$_{\Lambda}(U)$. So the latter assertion follows
from the former one.

(2) We proceed by using induction on $k$. The case $k=1$ is just
(1).

Now suppose $k \geq 2$. By induction assumption, we have that
$U$-lim.dim$_{\Lambda}(E_i) \leq i+1$  for any $0 \leq i \leq m-1$
and $U$-lim.dim$_{\Lambda}(E_{m+j}) \leq m+j$ for any $0 \leq j
\leq k-2$. By Corollaries 3.5 and 3.3, for any $N \in$mod $\Gamma
^{op}$, we have that s.grade$_U$Ext$_{\Gamma}^{i+1}(N, U) \geq i$
for any $1 \leq i \leq m$ and s.grade$_U$Ext$_{\Gamma}^{m+j}(N, U)
\geq m+j$ for any $1 \leq j \leq k-1$. By Corollary 3.3, it
suffices to prove that s.grade$_U$Ext$_{\Gamma}^{m+k}(N, U) \geq
m+k$.

Suppose that $N$ has a projective resolution as (8). By Lemma 3.4,
we have that $\mathcal{T}_U^i({\rm mod} \ \Gamma ^{op})=\Omega
_U^i({\rm mod} \ \Gamma ^{op})$ for any $1 \leq i \leq m+k$.
Notice that Coker$(Q_{m+k} \to Q_{m+k-1})$ is $(m+k-1)$-syzygy,
so, by Lemma 3.10, it is in $\Omega _U^{m+k-1}({\rm mod} \ \Gamma
^{op})$ and hence is in $\mathcal{T}_U^{m+k-1}({\rm mod} \ \Gamma
^{op})$, which implies that Ext$_{\Lambda}^i({\rm Tr}_U\Omega
_{\Gamma}^{m+k-1}(N), U)=0$ for any $1 \leq i \leq m+k-1$.

By assumption, $U$-lim.dim$_{\Gamma}(\bigoplus _{i=0}^mE'_i) \leq
m$ and $U$-lim.dim$_{\Gamma}(E'_{m+j}) \leq m+j$ for any $1 \leq j
\leq k-1$.  Then, by Corollary 3.3, we have that
s.grade$_U$Ext$_{\Lambda}^{m+k}(M, U) \geq m+k$ for any $M \in$mod
$\Lambda$.

Let $X$ be a submodule of Ext$_{\Gamma}^{m+k}(N, U)$. Then
grade$_UX \geq m+k-1$. By [9] Lemma 2, we have an exact sequence:
$$0 \to X \buildrel {f} \over \to {\rm Tr}_U\Omega
_{\Gamma}^{m+k-1}(N) \to {\rm Coker}f \to 0.$$ We then get an
embedding $0 \to {\rm Ext}_{\Lambda}^{m+k-1}(X, U) \to {\rm
Ext}_{\Lambda}^{m+k}({\rm Coker}f, U)$. Since
s.grade$_U$\linebreak ${\rm Ext}_{\Lambda}^{m+k}({\rm Coker}f, U)
\geq m+k$, grade$_U$Ext$_{\Lambda}^{m+k-1}(X, U) \geq m+k$. It
follows from Lemma 3.8 that grade$_UX \geq m+k$ and
s.grade$_U$Ext$_{\Gamma}^{m+k}(N, U) \geq m+k$. $\blacksquare$

\vspace{0.2cm}

Putting $m=0$, by Proposition 5.6(2), $U$-lim.dim$_{\Lambda}(E_i)
\leq i$ for any $0 \leq i \leq k-1$ if $U$-lim.dim$_{\Gamma}(E'_i)
\leq i$ for any $0 \leq i \leq k-1$. Combining this result with
Corollary 3.3(2) and their dual statements, we then get the
symmetry of $k$-Gorenstein modules (see [21] Theorem 7.5).

\vspace{0.2cm}

Putting $_{\Lambda}U_{\Gamma}={_{\Lambda}\Lambda _{\Lambda}}$, the
following corollary is an immediate consequence of Proposition
5.6, which is a generalization of the result of Auslander on the
symmetry of $k$-Gorenstein rings.

\vspace{0.2cm}

{\bf Corollary 5.7} {\it Assume that} $m$ {\it is a non-negative
integer and} {\it l.}fd$_{\Lambda}(I_i) \leq i+1$ {\it for any} $0
\leq i \leq m-1$.

(1) {\it If r.}fd$_{\Lambda}(\bigoplus _{i=0}^mI'_i) \leq m$, {\it
then l.}fd$_{\Lambda}(I_m) \leq m$; {\it Especially, if} {\it
l.}id$_{\Lambda}(\Lambda) \leq m$, {\it then l.}fd$_{\Lambda}(I_m)
\leq m$.

(2) {\it For a positive integer} $k$, {\it if
r.}fd$_{\Lambda}(\bigoplus _{i=0}^mI'_i) \leq m$ {\it and
r.}fd$_{\Lambda}(I'_{m+j}) \leq m+j$ {\it for any} $1 \leq j \leq
k-1$, {\it then l.}fd$_{\Lambda}(I_{m+j}) \leq m+j$ {\it for any}
$0 \leq j \leq k-1$.

\vspace{0.2cm}

When $m=0$, the result in Corollary 5.7(2) is equivalent to the
symmetry of $k$-Gorenstein rings (see [6] Theorem 3.7). In the
following, we give an example satisfying the conditions in
Corollary 5.7 for the case $m=1$ and $k=2$ as follows.

\vspace{0.2cm}

{\bf Example} Let $K$ be a field and $\Lambda$ a finite
dimensional $K$-algebra which given by the quiver:

$$
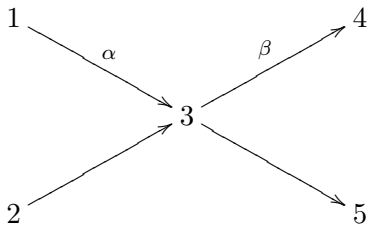
\xymatrix{1\ar[rrd]^\alpha &&&&4\\
&&3\ar[rru]^\beta\ar[rrd]\\
2\ar[rru]&&&&5}$$ modulo the ideal $\beta \alpha$. Then {\it
l.}fd$_{\Lambda}(I_0)$= {\it l.}fd$_{\Lambda}(I_1)$={\it
r.}fd$_{\Lambda}(I'_0)$={\it r.}fd$_{\Lambda}(I'_1)$=1, {\it
l.}fd$_{\Lambda}(I_2)$={\it r.}fd$_{\Lambda}(I'_2)$=2 and {\it
l.}id$_{\Lambda}(\Lambda)=${\it r.}id$_{\Lambda}(\Lambda)$=2.

\vspace{0.5cm}

{\bf Acknowledgements} The research of the author was partially
supported by Specialized Research Fund for the Doctoral Program of
Higher Education (Grant No. 20030284033). The author thanks Prof.
Takayoshi Wakamatsu for the useful discussion and suggestion. Some
part of this work was done during a visit of the author to Okayama
University from January to June, 2004. The author is grateful to
Prof. Yuji Yoshino for his kind hospitality.

\end{document}